\documentclass[12pt]{article}

\usepackage{amscd,amsmath,amssymb}
\usepackage[mathscr]{eucal}

\makeatletter \@addtoreset{equation}{section}\makeatother

\newtheorem{theorem}{Theorem}[section]
\newtheorem{lemma}[theorem]{Lemma}
\newtheorem{proposition}[theorem]{Proposition}

\title{\bf Berezin transform on the quantum unit ball}
\author{Dmitry Shklyarov$^\dagger$ \and Genkai Zhang$^\ddagger$}
\date{}

\begin{document}

\maketitle

\centerline{$^\dagger$\small{Institute for Low Temperature Physics \&
Engineering}}
\centerline {\small{47 Lenin Ave., 61103 Kharkov, Ukraine}}

\centerline{$^\ddagger$\small{Chalmers Tekniska H\"ogskola/G\"oteborgs
Universitet,
Matematik}} \centerline{\small{412 96, G\"oteborg, Sweden}}

\begin{center}
\small{e-mail: shklyarov@ilt.kharkov.ua, genkai@math.chalmers.se}
\end{center}

\medskip

\begin{center}\begin{minipage}[t]{4in}\small{ABSTRACT: We
introduce and study, in the framework of a theory of quantum Cartan domains,
a $q$-analogue of the Berezin transform on the unit ball. We construct
$q$-analogues of weighted Bergman spaces, Toeplitz operators and covariant
symbol calculus. In studying the analytical properties of the Berezin
transform we introduce also the $q$-analogue of the $SU(n,1)$-invariant
Laplace operator (the Laplace-Beltrami operator) and present related results
on harmonic analysis on the quantum ball.. These are applied to obtain an
analogue of one result by A.~Unterberger and H.~Upmeier. An explicit
asymptotic formula expressing the $q$-Berezin transform via the
$q$-Laplace-Beltrami operator is also derived. At the end of the paper, we
give an application of our results to basic hypergeometric $q$-orthogonal
polynomials.}
\end{minipage}
\end{center}

\medskip
\baselineskip 1.38pc

\section{Introduction}

\bigskip

Since appearance of quantum groups in the middle of 1980's there have been
different attempts to find an appropriate generalization to the $q$-case of
various classical constructions in analysis and geometry of Lie groups.
Recently, it became clear \cite{SV} that there should exist a substantial
$q$-analogue of the theory of Cartan domains (the most studied class of such
domains constitute the so-called classical domains \cite{Hua}). In turn,
this observation has opened a way to generalization of other important
theories about those domains. One of the most fascinating among them is the
Berezin's quantization \cite{B-cmp}. Though Berezin applied his construction
to a wide class of K\"ahler manifolds, the most complete and precise results
have being obtained just in the particular case of Cartan domains
\cite{B,UU}.

It should be noted that first attempts to find $q$-analogues of some
constructions of the Berezin's theory were made before the study of
$q$-Cartan domains has being initiated. For example, in \cite{KL} the
authors studied a two-parameter deformation $\mathcal{P}_{q,\lambda}$ of the
polynomial algebra $\mathbb{C}[z,\bar{z}]$ related to certain
$SU(1,1)$-covariant Poisson brackets on the unit disk. They showed among
other things that for any $\lambda$ the algebra $\mathcal{P}_{q,\lambda}$ is
acted upon by the quantum universal enveloping algebra $U_q\mathfrak{sl}_2$.
Also, for each $q$ the authors defined $q$-analogues of the weighted Bergman
spaces $L^2_a((1-|z|^2)^\lambda d\bar{z}dz)$ ($\lambda\geq0$) on the unit
disk and showed that $\mathcal{P}_{q,\lambda}$ could be realized as the
algebra of Toeplitz operators (with symbols from $\mathcal{P}_{q,0}$) on the
$q$-weighted Bergman space corresponding to $\lambda$. This observation
suggests that it is reasonable to regard $\mathcal{P}_{q,0}$ as a genuine
function algebra on the {\it quantum} unit disk and
$\mathcal{P}_{q,\lambda}$ as a result of Toeplitz quantization of
$\mathcal{P}_{q,0}$ which, in addition, respects the
$U_q\mathfrak{sl}_2$-actions in complete analogy with the classical setting.

This point of view was developed in details in \cite{bq1}. There, in
addition to the results of \cite{KL}, $q$-analogues of covariant symbols and
of the Berezin transform on the quantum unit disk was produced (see
\cite{bq} for a concise account). Also, in \cite{bq1} an explicit asymptotic
formula, expressing the $q$-Berezin transform via a $q$-Laplace-Beltrami
operator on the quantum disk, was derived. This allowed the authors to
produce an analogue of the Berezin's star product and thus to obtain a
formal deformation of the algebra $\mathcal{P}_{q,0}$.

The present paper continues the study of the $q$-Berezin transform initiated
in \cite{bq1}. Here we investigate the case of the quantum unit ball, whose
one-dimensional version is the aforementioned quantum unit disk, and
generalize to this case almost all constructions of \cite{bq1}. Namely, we
produce analogues of the weighted Bergman spaces, Toeplitz operators,
covariant symbols, and the Berezin transform on the unit ball. Also, we
define a $q$-analogue of the $SU(n,1)$-invariant Laplace operator (the
Laplace-Beltrami operator) and derive an explicit asymptotic formula
expressing the $q$-Berezin transform via the $q$-Laplace-Beltrami operator.

It is appropriate to mention here that in the case of the quantum unit ball
we encounter a new phenomenon imperceptible in the case of the quantum unit
disk. Namely, there is no any satisfactory analogue of coherent states on
the ball which are used to define covariant symbols of operators and
regarded as a basic ingredient of the Berezin's theory. An appropriate
analogue for the quantum unit disk was found in \cite{bq1}. However, that
was possible due to commutativity of the algebra of 'holomorphic functions'
on the quantum disk, but the commutativity fails in the case of the quantum
ball. Fortunately, even in the classical setting there is an alternative way
to define covariant symbols. Namely, the map $$
\begin{array}{ccl}\mathrm{operator}
\\
\mathrm{on\,\, a\,\, weighted\,\, Bergman\,\, space}
\end{array}\mapsto
\begin{array}{ccl}\mathrm{covariant\,\, symbol,}
\\
\mathrm{a\,\, function\,\, on\,\, the\,\, ball}
\end{array}
$$ is, roughly speaking, the adjoint of the map $$
\begin{array}{ccl}\mathrm{function}
\\
\mathrm{on\,\, the\,\, ball}
\end{array}\mapsto
\begin{array}{ccl}\mathrm{Toeplitz\,\, operator}
\\
\mathrm{on\,\, a\,\, weighted\,\, Bergman\,\, space}
\end{array}
$$ with respect to certain $SU(n,1)$-invariant inner products in the spaces
of functions and operators (see \cite{UU}). The significance of the Toeplitz
and covariant calculi is of course well known and has been intensively
studied. We exploit just this idea to define covariant symbols and thus the
Berezin transform in the $q$-setting.

Let us turn now to description of the contents of the paper.

Our results rely heavily on function theory in the quantum unit ball. In
Section \ref{2} we recall the basic setup and results from that theory
obtained earlier by L.~Vaksman and his group \cite{SSV1, SSV2} in the more
general setting of quantum {\it matrix} balls. In brief, we define an
involutive algebra of 'polynomials' on the quantum $\mathbb{C}^n$ which
generalize the algebra $\mathcal{P}_{q,0}$ mentioned above. The algebra is
endowed with an action of the quantum universal enveloping algebra
$U_q\mathfrak{su}_{n,1}$. This is a counterpart of the classical
$SU(n,1)$-action on the unit ball. Using the polynomial algebra, we produce
the spaces of 'finite functions' and 'distributions' on the quantum ball
which inherit the $U_q\mathfrak{su}_{n,1}$-action. The crucial property,
which justifies our definition of the space of finite functions on the
quantum ball, is the existence of a $U_q\mathfrak{su}_{n,1}$-invariant
integral on that space, an analogue of the integral with respect to the
$SU(n,1)$-invariant measure on the ball. We present an explicit formula for
that integral. We note that the existence of the integral allows one to
formulate many problems of harmonic analysis on the quantum ball. Though the
present paper does not deal with harmonic analysis, we develop
 some partial results (see description of
Section \ref{4} below) which are needed to put
our result in a coherent and rigorous mathematics
context.

Section \ref{3} is devoted to
 find  an appropriate algebraic
and analytical setting.
 There we start developing of the
$q$-Berezin's theory and present all necessary definitions. First, we define
analogues of the weighted Bergman spaces in the ball and produce the notion
of Toeplitz operators with polynomial or finite symbols. Though the finite
symbols are of more importance for our exposition, we devote a separate
subsection to the study of the algebra of Toeplitz operators with polynomial
symbols. This is motivated by the observation that this algebra is a
multivariable generalization of the algebra $\mathcal{P}_{q,\lambda}$
studied in \cite{KL}. Results of \cite{KL} suggest that the algebra of
Toeplitz operators with polynomial symbols might be interesting from the
point of view of the theory of operator algebras. In the last subsection we
present definitions of covariant symbols and of the $q$-Berezin transform.
Thus, in Section \ref{3} we define three important maps: the map that sends
a function to the corresponding Toeplitz operator, the map that sends an
operator to its covariant symbol (a function on the quantum ball), and,
finally, their composition - the $q$-Berezin transform - which sends a
function to the covariant symbol of the corresponding Toeplitz operator. We
observe also that the three maps intertwine natural
$U_q\mathfrak{su}_{n,1}$-actions on the spaces of functions and operators.
This quantum group symmetry of the $q$-Berezin transform simplifies
drastically many computations in subsequent sections.

To investigate the $q$-Berezin transform further, for instance, to find its
asymptotic expansion, we need some results on harmonic analysis,
particularly, on $q$-spherical transform in the quantum ball. These are
presented in Section \ref{4}. First of all, we define there an analogue of
the Laplace-Beltrami operator on the ball. At first sight, our definition
seems to be somewhat unusual: the $q$-Laplace-Beltrami operator appears as
the first term in the asymptotic expansion of the $q$-Berezin transform.
However, this agrees with the classical Berezin's theory. Besides, we
believe that results of Section \ref{4} provide a sufficient evidence of
rationality of such a definition. In particular, a 'radial part' of our
$q$-Laplace-Beltrami operator turns out to be a second order difference
operator which tends to the radial part of the classical Laplace-Beltrami
operator as $q$ approaches 1. Moreover, the difference operator turned out
to be quite well studied by experts in the theory of $q$-special functions
(see, for instance, \cite{KoeSt}). In particular, the problem of expansion
in eigenfunctions of that difference operator had being already solved, and
we recall explicit formulas in the text. In view of our approach, those
eigenfunctions should be regarded as analogues of the spherical functions on
the ball while the expansion is an analogue of the spherical transform (see
\cite{Helg}). Moreover, the Plancherel measure, which appears in the inverse
transform, involves an analogue of the Harish-Chandra $c$-function for the
ball.

After we present the results on the $q$-spherical transform, we turn back in
Section \ref{5} to the study of the $q$-Berezin transform. In the first
subsection of Section \ref{5} we consider its 'radial part'. We prove that
the radial part of the transform is extended to a bounded selfadjoint
operator on the space of square-integrable radial function on the quantum
ball (we mean square-integrability with respect to the invariant
$q$-integral). Moreover, the resulting operator commutes with the radial
part of the $q$-Laplace-Beltrami operator which is also a bounded
selfadjoint operator on the same space. Since the latter operator has simple
spectrum (this is stated in Section \ref{4}), the radial part of the
$q$-Berezin transform appears to be a function of the radial part of the
$q$-Laplace-Beltrami operator. We found the function explicitly. We note
that in the classical case the corresponding result was obtained in
\cite{UU} for any Cartan domain. This computation has a number of
consequences some of which we use further in Section \ref{6}. In the last
subsection of Section \ref{5}, we present an asymptotic formula expressing
the $q$-Berezin transform via the $q$-Laplace-Beltrami operator.

Section \ref{6}, the last section of the paper, is meant for those readers
who are interested in various mathematical applications of the Berezin
transform. It is however nevertheless related to quantization. In the theory
of Wick quantization one associates to each function an operator so that the
one-dimensional projections onto the coherent state at a point (in the phase
space) is assigned to the delta function. However if the given Hilbert space
has on orthogonal basis the projection onto the basis vectors are also
natural subjects so it is interesting to compute their covariant symbols.
Here in our setting we show that the spherical transforms of those covariant
symbols are roughly speaking hypergeometric orthogonal $q$-polynomials. We
derive orthogonality relations for those orthogonal polynomials, the
so-called continuous dual $q$-Hahn polynomials \cite{KS}. The entire family
of these polynomials depends on three parameters (not counting $q$) while we
treat a two-parameter sub-family which appears naturally in connection with
the $q$-Berezin transform. The idea we use appeared for the first time in
\cite{jp+gkz-pl3} in the classical setting (see also \cite{Ne}).

\medskip

Let us make some comments about the most important agreements and notations
used in the paper. First of all, the quantum-group parameter $q$ is supposed
to be a number from the interval $(0\,,\,1)$. Next, we put $h:=\log q^{-2}$,
so $h>0$; and $t:=q^{2\alpha}$ will be the deformation parameter in the
Berezin transform. Finally, we describe some convenient multi-index
notations. We denote multi-indices by underlined small letters, the zero
multi-index by $\underline{0}$. For $\underline{i}=(i_1,\ldots, i_n)$,
$\underline{j}=(j_1,\ldots, j_n)$ we put
$|\underline{i}|=i_1+i_2+\ldots+i_n$ and $
\underline{i}\times\underline{j}=(i_1\cdot j_1,\ldots,i_n\cdot j_n).$ For
non-commuting variables $z_1,z_2,\ldots, z_n$ and $z^*_1,z^*_2,\ldots,
z^*_n$ we set $$ \mathbf{z}^{\underline{i}}:=z_n^{i_n}\cdot
z_{n-1}^{i_{n-1}}\cdot\ldots\cdot z_1^{i_1}, \qquad
\mathbf{z}^{*\underline{i}}:=z_1^{*i_1}\cdot\ldots\cdot
z_{n-1}^{*i_{n-1}}\cdot z_n^{*i_n}.$$ Finally, we denote
$P(n)=\{\underline{k}\in\mathbb{Z}^n_{\geq0}\;|\;k_1\geq k_2\geq\ldots\geq
k_n\}. $

\medskip

{\bf Acknowledgments.} It is our pleasure to thank Leonid Vaksman for many
interesting discussions. This research was supported by Royal Swedish
Academy of Sciences under the program 'Cooperation between researchers in
Sweden and the former Soviet Union'. We thank also the anonymous referee for
expert comments and advices.

\bigskip

\section{Basics of function theory in the quantum unit ball}\label{2}

The aim of Section \ref{2} is to introduce some notions of function theory
in the quantum ball. We define the involutive algebra of finite functions
and the space of distributions on the $q$-ball. These spaces admit actions
of certain quantum universal enveloping algebra. Finally, we produce an
explicit formula for the invariant integral on the space of finite
functions. The material presented in this section is not new: the quantum
ball treated here is a particular case of the quantum {\sl matrix} ball
considered in details in \cite{SSV1,SSV2}. Most results can be found in
those papers, however, we reformulate some of them for our particular
purpose.

\medskip

\subsection{Polynomials on the quantum ${\mathbb C}^n$}

An initial object in constructing function theory in the quantum ball is the
unital involutive algebra given by its generators $z_1, z_2,\ldots z_n$ and
the relations

\begin{equation}
\label{zz} z_iz_j=qz_jz_i, \quad i<j
\end{equation}
\begin{equation}\label{z*z} z_i^* z_j=qz_jz_i^*, \quad i\ne j
\end{equation}
\begin{equation} \label{z*z1}
z_j^*z_j=q^2z_jz_j^*+(1-q^{2})(1-\sum_{k=j+1}^nz_kz_k^*).
\end{equation}
We denote this algebra by $\mathcal{P}(\mathbb{C}^n)_q$ and treat it as the
polynomial algebra on the quantum $\mathbb{C}^n$. The algebra
$\mathcal{P}(\mathbb{C}^n)_q$ is a particular case of the algebra
$\mathrm{Pol}(\mathrm{Mat}_{m,n})_q$ of polynomials on the quantum space of
$m\times n$ matrices considered in \cite {SSV1, SSV2}. Specifically,
$\mathcal{P}(\mathbb{C}^n)_q$ coincides with
$\mathrm{Pol}(\mathrm{Mat}_{1,n})_q$. (Moreover, it is not hard to show that
$\mathcal{P}(\mathbb{C}^n)_q$ is isomorphic to the known twisted CCR-algebra
\cite{PW}.)

In the sequel we use the algebra $\mathcal{P}(\mathbb{C}^n)_q$ to produce
some 'functional' spaces and algebras on the quantum ball. Those spaces, as
well as the algebra $\mathcal{P}(\mathbb{C}^n)_q$ itself, are endowed with
structures of $U_q\mathfrak{su}_{n,1}$-modules where
$U_q\mathfrak{su}_{n,1}$ is certain $*$-Hopf algebra called the quantized
universal enveloping algebra for $\mathfrak{su}_{n,1}$. Let us turn to
precise definitions.

We recall first the definition of $U_q\mathfrak{su}_{n,1}$. It is a 'real
form' of the Drinfeld-Jimbo quantum universal enveloping algebra
$U_q\mathfrak{sl}_{n+1}$. The latter algebra is the unital algebra given by
the generators $E_i, F_i, K_i^{\pm1}$, $i=1,2,\ldots n$, and the following
relations

$$
K_iK_j=K_jK_i,\qquad K_iK_i^{-1}=K_i^{-1}K_i=1,\qquad
K_iE_j=q^{a_{ij}}E_jK_i,
$$ $$ K_iF_j=q^{-a_{ij}}F_jK_i, \qquad
E_iF_j-F_jE_i=\delta_{ij}(K_i-K_i^{-1})/(q-q^{-1}) $$ $$
E_i^2E_j-(q+q^{-1})E_iE_jE_i+E_jE_i^2=0,\qquad |i-j|=1 $$ $$
F_i^2F_j-(q+q^{-1})F_iF_jF_i+F_jF_i^2=0,\qquad |i-j|=1 $$ $$
[E_i,E_j]=[F_i,F_j]=0,\qquad |i-j|\ne 1 $$ with $(a_{ij})$ being the Cartan
matrix of type $A_n$. Moreover, $U_q\mathfrak{sl}_{n+1}$ is a Hopf algebra.
The comultiplication $\Delta$, the antipode $S$, and the counit
$\varepsilon$ are determined by $$ \Delta(E_i)=E_i \otimes 1+K_i \otimes
E_i,\quad \Delta(F_i)=F_i \otimes K_i^{-1}+1 \otimes F_i,\quad
\Delta(K_i)=K_i \otimes K_i, $$ $$ S(E_i)=-K_i^{-1}E_i,\qquad
S(F_i)=-F_iK_i,\qquad S(K_i)=K_i^{-1}, $$
$$\varepsilon(E_i)=\varepsilon(F_i)=0,\qquad \varepsilon(K_i)=1.$$

The quantum universal enveloping algebra $U_q\mathfrak{su}_{n,1}$ is defined
as the $*$-Hopf algebra $(U_q\mathfrak{sl}_{n+1},*)$ with $*$ being the
involution in $U_q\mathfrak{sl}_{n+1}$ given by
\begin{align*}
E_n^*&=-K_nF_n,&F_n^*&=-E_nK_n^{-1},&(K_n^{\pm 1})^*&=K_n^{\pm
1},&
\\ E_j^*&=K_jF_j,&F_j^*&=E_jK_j^{-1},&(K_j^{\pm 1})^*&=K_j^{\pm1},
&\text{for}\quad j \ne n.
\end{align*}
(see \cite{ChP} for basic definitions concerning $*$-Hopf algebras).

We need also the notion of $U_q\mathfrak{su}_{n,1}$-module algebra. Let $A$
be a Hopf algebra. An algebra $F$ is said to be an $A$-module algebra if $F$
carries a structure of $A$-module and multiplication in $F$ agrees with the
$A$-action, i.e. the multiplication $F\otimes F\rightarrow F$ is a morphism
of $A$-modules.\footnote{We recall that for any $A$-modules $V_1$, $V_2$
their tensor product is endowed with an $A$-module structure via the
comultiplication $\Delta:A\to A\otimes A$.} If $A$ or $F$ has some extra
structures this definition includes natural additional requirements. For
example, in the case of a unital algebra $F$ one adds the requirement of
$A$-invariance of the unit: $\xi(1)=\varepsilon(\xi)\cdot1$, $\xi\in A$. In
the case of an involutive algebra $F$ and a $*$-Hopf algebra $A$ one imposes
the requirement of agreement of the involutions:
\begin{equation}\label{agree}
(\xi(f))^*=S(\xi)^*(f^*),\quad\xi\in A, f\in F.
\end{equation}

Some natural examples of module algebras appear in the classical setting.
Suppose $X$ is a smooth $G$-space with $G$ being a Lie group. The $G$-action
induces an action of the universal enveloping algebra $U\mathfrak{g}$ in the
space $C^\infty(X)$ via differential operators. The usual Leibnitz rule for
the differentiation of product of two functions means that $C^\infty(X)$ is
a $U\mathfrak{g}$-module algebra. This example suggests the use of the
language of module algebras for description of algebras of functions on {\sl
quantum} $G$-spaces.

The structure of $U_q\mathfrak{su}_{n,1}$-module algebra in
$\mathcal{P}(\mathbb{C}^n)_q$ which we are going to present has the
following classical counterpart. The unit ball
$U_n=\{\mathbf{z}\in\mathbb{C}^n | \|\mathbf{z}\|<1\}$ is a homogeneous
space of the group $SU(n,1)$ whose elements act via (biholomorphic)
linear-fractional transformations. Elements of the universal enveloping
algebra $U\mathfrak{su}_{n,1}$ act on the space $C^\infty(U_n)$ via
differential operators with polynomial coefficients and thus keep the
subspace of polynomials invariant. This induces a
$U\mathfrak{su}_{n,1}$-module algebra structure in the polynomial algebra on
$\mathbb{C}^n$. We turn now to the quantum case.

To describe the $U_q\mathfrak{su}_{n,1}$-module algebra structure in
$\mathcal{P}(\mathbb{C}^n)_q$ we consider the action of the elements $E_i,
F_i, K_i^{\pm1}$ on the generators $z_i, z_i^*$. Moreover, by (\ref{agree})
it is sufficient to present formulas for the action of $E_i, F_i,
K_i^{\pm1}$ on the 'holomorphic' part $z_1,\ldots z_n$ of generators of
$\mathcal{P}(\mathbb{C}^n)_q$. These are given in the following

\medskip

\begin{proposition}\label{2.1} There exists a unique structure of
$U_q\mathfrak{su}_{n,1}$-module algebra in $\mathcal{P}(\mathbb{C}^n)_q$
such that, for $k\ne n$,
\begin{equation}\label{vak-h's}
K_kz_i=\left \{\begin{array}{ccl}qz_i, &i=k
\\
q^{-1}z_i, & i=k+1
 \\z_i, \quad  & {\rm otherwise,}\end{array}\right.
\end{equation}
\begin{equation}\label{vak-f's}
F_kz_i=q^{1/2}\cdot \left
\{\begin{array}{ccl} z_{i+1}, & i=k
 \\ 0, \quad &{\rm
otherwise,}\end{array}\right.
\end{equation}
\begin{equation}\label{vak-e's}
E_kz_i=q^{-1/2}\cdot \left \{\begin{array}{ccl}z_{i-1}, &  i=k+1
\\ 0,& {\rm otherwise,}\end{array}\right.
\end{equation}
and
\begin{equation}\label{vak-h}
K_nz_i=\left \{\begin{array}{ccl}q^2z_i, &i=n
\\ qz_i, &{\rm otherwise,}\end{array}\right.
\end{equation}
\begin{equation}\label{vak-f}
F_nz_i=q^{1/2}\cdot \left \{\begin{array}{ccl}1, & i=n  \\ 0
, &{\rm otherwise,}\end{array}\right.
\end{equation}
\begin{equation}\label{vak-e}
E_nz_i=-q^{1/2}\cdot \left \{\begin{array}{ccl} z_n^2, & i=n
 \\ z_nz_i, &{\rm
otherwise.}\end{array}\right.
\end{equation}

\end{proposition}

\medskip

\noindent This statement may be deduced from Proposition 2.1 and Corollary
5.6 in \cite{SSV1} by easy computations. In the next subsection we will
introduce some other important $U_q\mathfrak{su}_{n,1}$-module algebras.

\medskip

\subsection{Finite functions on the quantum ball}\label{2.2}

It should be noted that in the classical case the space of polynomials does
not suit purposes of harmonic analysis in the ball because the volume of the
ball with respect to the $SU(n,1)$-invariant measure
\begin{equation}\label{measure}
d\nu(\mathbf{z})=\frac{dm(\mathbf{z})}{(1-\|\mathbf{z}\|^2)^{n+1}}
\end{equation} (where $dm(\mathbf{z})$ is the normalized Lebesgue measure) is
infinite. One observes the same problem in the quantum setting: there is no
$U_q\mathfrak{su}_{n,1}$-invariant integral on
$\mathcal{P}(\mathbb{C}^n)_q$, i. e. a linear positive functional
$\nu:\mathcal{P}(\mathbb{C}^n)_q\rightarrow\mathbb{C}$ such that $
\nu(\xi(f))=\varepsilon(\xi)\cdot \nu(f)$ for any $\xi\in
U_q\mathfrak{su}_{n,1}$ ($\varepsilon$ is the counit of
$U_q\mathfrak{su}_{n,1}$). Thus it is reasonable to look for a quantum
counterpart of the space of finite functions on the ball (functions with
compact supports inside the ball). The following construction for such a
counterpart was proposed in \cite[Section 7]{SSV1}.

Let us add to the algebra $\mathcal{P}(\mathbb{C}^n)_q$ one more generator
$f_0$ which satisfies the relations
\begin{equation}\label{f0}
f_0^*=f_0^2=f_0,\qquad z_i^*f_0=f_0z_i=0,\quad i=1,2,\ldots n.
\end{equation}
The resulting involutive algebra will be denoted by $\mathcal{F}(U_n)_q$. It
is demonstrated in \cite[Section 7]{SSV1} that one may extend the structure
of $U_q\mathfrak{su}_{n,1}$-module algebra in $\mathcal{P}(\mathbb{C}^n)_q$
to one in $\mathcal{F}(U_n)_q$ as follows
\begin{equation}\label{ef}
K_kf_0=f_0, \qquad F_kf_0=E_kf_0=0
\end{equation}
with $k \ne n$ and
\begin{equation}\label{hfe}
K_nf_0=f_0,\qquad F_nf_0=-\frac{q^{1/2}}{q^{-2}-1}f_0 \cdot z_n^*,\qquad
E_nf_0=-\frac{q^{1/2}}{1-q^2}z_n \cdot f_0.
\end{equation}
The involutive algebra $\mathcal{D}(U_n)_q$ of finite functions on the
quantum ball is defined as the two-sided ideal $\mathcal{F}(U_n)_q\cdot
f_0\cdot \mathcal{F}(U_n)_q$ in $\mathcal{F}(U_n)_q$. Due to (\ref{ef}),
(\ref{hfe}), $\mathcal{D}(U_n)_q$ is a $U_q\mathfrak{su}_{n,1}$-module
subalgebra in $\mathcal{F}(U_n)_q$. In the next subsection we present an
explicit formula for a $U_q\mathfrak{su}_{n,1}$-invariant integral on
$\mathcal{D}(U_n)_q$.

The above definition of $\mathcal{D}(U_n)_q$ is convenient for many purposes
but performing computations. We present therefore an alternative description
of $\mathcal{D}(U_n)_q$. For that purpose we let $H$ be the
$\mathcal{P}(\mathbb{C}^n)_q$-module given by its unique generator $e_{0}$
such that $z_i^*e_{0}=0$, $i=1,2,\ldots n$. By the formulas (\ref{z*z}) and
(\ref{z*z1}), $H=\mathcal{P}(\mathbb{C}^n)_q\cdot e_{0}=
\mathbb{C}[\mathbb{C}^n]_q\cdot e_{0}$ with $\mathbb{C}[\mathbb{C}^n]_q$
being the unital (non-involutive) subalgebra in
$\mathcal{P}(\mathbb{C}^n)_q$ generated by $z_i$, $i=1,2,\ldots n$.
Moreover, the elements $\{\mathbf{z}^{\underline{m}}\cdot
e_{0}\}_{\underline{m}\in\mathbb{Z}^n_{\geq0}}$ constitute a basis in $H$.
The following statement may be deduced from known results concerning the
twisted CCR-algebra (see, for instance, \cite{Prosk}).

\medskip

\begin{proposition}
 i) There exists a (unique up to a positive multiplier) scalar
product in $H$ such that $$ (f\cdot e_1, e_2)=(e_1, f^*\cdot e_2)$$ for any
$f\in\mathcal{P}(\mathbb{C}^n)_q$ and $e_1,e_2\in H$.

ii)The corresponding $*$-representation $T$ of $\mathcal{P}(\mathbb{C}^n)_q$
in the completion $\overline{H}$ of the pre-Hilbert space $H$ (the so-called
Fock representation) is a faithful irreducible representation by bounded
operators.
\end{proposition}

\medskip

We let $C^*(\mathcal{P}(\mathbb{C}^n)_q)$ be the $C^*$-algebra generated by
$\mathcal{P}(\mathbb{C}^n)_q$ via the representation $T$. We describe now an
embedding of $\mathcal{D}(U_n)_q$ in $C^*(\mathcal{P}(\mathbb{C}^n)_q)$.

The representation $T$ can be extended to a faithful $*$-representation of
the algebra $\mathcal{F}(U_n)_q$ by setting $$ T(f_0)=\mathrm{orthogonal \;
projection\; onto}\quad\mathbb{C}\cdot e_{0}.$$ It is easy to show that the
subalgebra $\{T(f)\;|\;f\in\mathcal{D}(U_n)_q\}$ of the algebra of bounded
operators on $\overline{H}$ coincides with the subalgebra of those operators
whose matrices in the basis $\{\mathbf{z}^{\underline{m}}\cdot
e_{0}\}_{\underline{m}\in\mathbb{Z}^n_{\geq0}}$ have only finitely many
non-zero entries. Thus the algebra $\mathcal{D}(U_n)_q$ can be realized as a
subalgebra of $C^*(\mathcal{P}(\mathbb{C}^n)_q)$. It admits the following
transparent description.

Let $y_i{=}1-z_iz_i^*-\ldots -z_nz_n^*\in\mathcal{P}(\mathbb{C}^n)_q,$
$i=1,2,\ldots n.$ Then $(y_1, \ldots, y_n)$ is a tuple of pair-wise
commuting positive operators on $\overline{H}$. The joint spectrum of the
tuple $(y_1,,\ldots, y_n)$ in $C^*(\mathcal{P}(\mathbb{C}^n)_q)$ is the
closure in $\mathbb{R}^n$ of the '$q$-simplex'
\begin{equation}\label{simplex} \mathfrak{M}=\{(q^{2k_1}, q^{2k_2},\ldots
q^{2k_n})\;|\; \underline{k}\in P(n)\}\end{equation} ($P(n)$ is defined in
the Introduction). This is a consequence of the commutation relations
\begin{equation}\label{yi}
z_iy_j=\left \{\begin{array}{ccl} q^{-2}y_jz_i, & i\le j
 \\ y_jz_i, &{\rm
otherwise,}\end{array}\right. \qquad z^*_iy_j=\left \{\begin{array}{ccl}
q^{2}y_jz^*_i, & i\le j
 \\ y_jz^*_i, &{\rm
otherwise.}\end{array}\right.
\end{equation}
We will hence forth identify a function $f(y_1,\dots, y_n)$ with a function
$f$ on the set $\mathfrak{M}$ via the spectral calculus.

Using the definition of $y_i$'s and (\ref{yi}), one can write an arbitrary
element $f\in\mathcal{P}(\mathbb{C}^n)_q$ uniquely in the form
\begin{equation}\label{form}
f=\sum_{\underline{i}\times\underline{j}=0}
\mathbf{z}^{\underline{i}}f_{\underline{i},\underline{j}}(y_1,y_2,\ldots
y_n)\mathbf{z}^{*\underline{j}}.
\end{equation}
The subalgebra $\mathcal{D}(U_n)_q\subset C^*(\mathcal{P}(\mathbb{C}^n)_q)$
may be identified with the algebra of finite sums of the form (\ref{form})
whose coefficients are functions on $\mathfrak{M}$ with finite supports.

We will use the both descriptions of $\mathcal{D}(U_n)_q$. It is clear that
the 'distinguished' finite function $f_0$ which appears in the definition of
$\mathcal{D}(U_n)_q$ may be described, in those terms, as follows: $
f_{\underline{i},\underline{j}}\equiv0$ provided
$\underline{i}\neq\underline{0}$ or $\underline{j}\neq\underline{0}$, and

$$ f_{\underline{0},\underline{0}}(q^{2k_1}, q^{2k_2},\ldots q^{2k_n})=
\left \{\begin{array}{ccl}1, &
\underline{k}=\underline{0}
\\ 0, & {\rm otherwise.}\end{array}\right.$$

The following crucial property of $f_0$ will simplify proofs of many results.

\medskip

\begin{proposition}\label{f00}
$f_0$ generates $\mathcal{D}(U_n)_q$ as a $U_q\mathfrak{su}_{n,1}$-module.
\end{proposition}

\medskip

\noindent A proof can be found in \cite[Section 7]{SSV1}.

\medskip

\subsection{Invariant integral} In this subsection we present an explicit
formula for the $U_q\mathfrak{su}_{n,1}$-invariant integral on the quantum
ball. This formula was found in \cite[Section 9]{SSV1}.

Remind the notation $T$ for the Fock representation of
$\mathcal{P}(\mathbb{C}^n)_q$ and $H$ for a
$\mathcal{P}(\mathbb{C}^n)_q$-module, a dense linear subspace in the Hilbert
space $\overline{H}$ of th representation $T$ (see the previous subsection).
It is easy to observe that $T(f)(H)\subset H$ for any finite function $f$.
Thus $H$ become a $\mathcal{D}(U_n)_q$-module. It is convenient to identify
the $\mathcal{D}(U_n)_q$-module $H$ with the left ideal
$\mathbb{C}[\mathbb{C}^n]_q\cdot f_0$ (this is possible due to (\ref{f0})),
with which we may equip $H$ with some extra structures. In particular, the
isomorphism of $\mathcal{D}(U_n)_q$-modules
$$H\simeq\mathbb{C}[\mathbb{C}^n]_q\cdot f_0\subset\mathcal{D}(U_n)_q$$ and
formulas (\ref{vak-h's}), (\ref{vak-h}), (\ref{ef}), (\ref{hfe}) define $H$
as a $U_q\mathfrak{h}$-module where $U_q\mathfrak{h}$ is the Hopf subalgebra
in $U_q\mathfrak{sl}_{n+1}$ generated by $K_i^{\pm1}$, $i=1,2,\ldots n$.
Denote the corresponding representation of $U_q\mathfrak{h}$ in $H$ by
$\Gamma$. The following statement was proved in \cite{SSV1}.

\medskip

\begin{proposition}\label{integral}
The linear functional on $\mathcal{D}(U_n)_q$ given by $$ f\mapsto \nu(f){=}
\mathrm{tr}\left(T(f)\Gamma\left(\prod_{j=1}^nK_j^{-j(n+1-j)}\right)
\right)$$ is a $U_q\mathfrak{su}_{n,1}$-invariant integral, i.e.
$\nu(f^*f)>0$, $f\ne0$, and $\nu(\xi(f))=\varepsilon(\xi)\cdot \nu(f)$,
$\xi\in U_q\mathfrak{su}_{n,1}$.
\end{proposition}

\medskip

\noindent We shall use the normalized integral
$$\int_{U_n}fd\nu_q=(q^2;q^2)_n\cdot\nu(f)$$ with
$(q^2;q^2)_n{=}(1-q^2)(1-q^4)\ldots(1-q^{2n})$.

The invariant integral is unique up to a positive scalar. This is an
immediate consequence of Proposition \ref{f00}.

We can rewrite now the formula for the invariant integral in a more
convenient form. Suppose $f\in\mathcal{D}(U_n)_q$ is written in the form
(\ref{form}). Then one checks easily that
\begin{equation}\label{invint}
\int_{U_n}fd\nu_q=(q^2;q^2)_n\cdot\sum_{\underline{k}\in P(n)
}f_{\underline{0},\underline{0}}(q^{2k_1}, q^{2k_2},\ldots q^{2k_n})\cdot
q^{-2nk_1}\cdot q^{2k_2}\cdot\ldots\cdot q^{2k_n}.
\end{equation}

In the next proposition we explain how the
$U_q\mathfrak{su}_{n,1}$-invariant integral makes $\mathcal{D}(U_n)_q$ into
a unitary $U_q\mathfrak{su}_{n,1}$-module and thus provides a setup for
harmonic analysis in the quantum ball.

\medskip

\begin{proposition}\label{scpro}
The scalar product in $\mathcal{D}(U_n)_q$ given by $$(\varphi_1,
\varphi_2)=\int_{U_n}\varphi^*_2\cdot \varphi_1d\nu_q$$ is
$U_q\mathfrak{su}_{n,1}$-invariant, i.e.
$(\xi(\varphi_1),\varphi_2)=(\varphi_1,\xi^*(\varphi_2))$ for any
$\varphi_1,\varphi_2\in \mathcal{D}(U_n)_q$ and any $\xi\in
U_q\mathfrak{su}_{n,1}.$
\end{proposition}

\medskip

\noindent{\bf Proof.} This follows from the condition (\ref{agree}) in
$\mathcal{D}(U_n)_q$ and the following formula of 'integrating by parts':
for any $\varphi_1,\varphi_2\in\mathcal{D}(U_n)_q$ and $\xi\in
U_q\mathfrak{su}_{n,1}$ one has

\begin{equation}\label{inparts}
\int_{U_n}\xi(\varphi_1)\cdot\varphi_2d\nu_q=\int_{U_n}\varphi_1\cdot
S(\xi)(\varphi_2)d\nu_q
\end{equation}
with $S$ being the antipode of $U_q\mathfrak{su}_{n,1}$. To prove the latter
formula, it is sufficient to verify the equality on generators $E_i, F_i,
K_i^{\pm1}$ of $U_q\mathfrak{su}_{n,1}$. In this particular case the
equality is equivalent to the $U_q\mathfrak{su}_{n,1}$-invariance of the
integral. \hfill $\blacksquare$

\medskip

\subsection{Distributions on the quantum ball}\label{2.4}

The aim of this subsection is to define the space $\mathcal{D}(U_n)'_q$ of
distributions on the quantum ball. We follow ideas in \cite[Section 1]{qd1}
where the simplest case $n=1$ is treated.

Equip the vector space $\mathcal{D}(U_n)_q$ with the weakest topology so
that all the linear functionals $$
l_{\underline{i},\underline{j}}^{\underline{k}}:f\mapsto
f_{\underline{i},\underline{j}}(q^{2k_1}, q^{2k_2},\ldots q^{2k_n})$$ are
continuous (here $\underline{i},\underline{j}\in\mathbb{Z}_{\geq0}^n$
satisfy $\underline{i}\times\underline{j}=\underline{0}$, $\underline{k}\in
P(n)$, and $f_{\underline{i},\underline{j}}$ is the corresponding
coefficient in the expansion (\ref{form}) of $f$).

The following results can be proved by direct however tedious computations,
which we omit.

\medskip

\begin{proposition}
For a fixed $f\in\mathcal{D}(U_n)_q$ the operator in $\mathcal{D}(U_n)_q$ of
the right multiplication by $f$ is continuous. Also, the
$U_q\mathfrak{su}_{n,1}$-action in $\mathcal{D}(U_n)_q$ is continuous.
\end{proposition}

\medskip

Let $\mathcal{D}(U_n)'_q$ be the completion of the topological vector space
$\mathcal{D}(U_n)_q$. We shall use the concrete realization of
$\mathcal{D}(U_n)'_q$ as the space of infinite sums of the form (\ref{form})
whose coefficients $f_{\underline{i},\underline{j}}$ are arbitrary functions
on the $q$-simplex $\mathfrak{M}$ (see (\ref{simplex})). The above
proposition, by continuity, implies that $\mathcal{D}(U_n)'_q$ is a right
$\mathcal{D}(U_n)_q$-module and a $U_q\mathfrak{su}_{n,1}$-module. Moreover,
these structures agree: the multiplication map
$\mathcal{D}(U_n)'_q\otimes\mathcal{D}(U_n)_q\rightarrow\mathcal{D}(U_n)'_q$
is a morphism of $U_q\mathfrak{su}_{n,1}$-modules. This is a consequence of
$U_q\mathfrak{su}_{n,1}$-moduleness of the algebra $\mathcal{D}(U_n)_q$.

Let us construct a non-degenerate pairing
$\mathcal{D}(U_n)'_q\times\mathcal{D}(U_n)_q\rightarrow\mathbb{C}$ which
will justify the term 'distribution'. Let $\{\psi_\beta\}_{\beta}$ be a
convergent net in $\mathcal{D}(U_n)_q$ and suppose
$\lim_{\beta}\psi_\beta=\psi$, $\psi\in\mathcal{D}(U_n)'_q$. Then for any
$\varphi\in\mathcal{D}(U_n)_q$ the limit $\lim_{\beta}\int_{U_n}\psi_\beta
\varphi d\nu_q$ exists and depends only on $\psi$ and $\varphi$ (i. e.
independent of a choice of $\psi_\beta$). Indeed, let us prove existence of
$\lim\int_{U_n}\psi_\beta \varphi d\nu_q$. Proposition \ref{f00} implies
that there exists $\xi_\varphi\in U_q\mathfrak{su}_{n,1}$ such that
$\varphi=\xi_\varphi(f_0)$. By (\ref{inparts}) $$ \int_{U_n}\psi_\alpha
\varphi d\nu_q=\int_{U_n}S^{-1}(\xi_\varphi)(\psi_\beta) f_0 d\nu_q.$$ Since
the $U_q\mathfrak{su}_{n,1}$-action is continuous, it is sufficient now to
prove existence of the limit in the particular case $\varphi:=f_0$. But in
this case existence equivalent to continuity of the functional
$l_{\underline{0}\underline{0}}^{\underline{0}}$ and thus is obvious.
Independence of the limit on a particular choice of the net $\{\psi_\beta\}$
can be proved by the same arguments.

Thus we get a well defined pairing
$$\mathcal{D}(U_n)'_q\times\mathcal{D}(U_n)_q\rightarrow\mathbb{C},\qquad
(\psi,\varphi)\mapsto \int_{U_n}\psi \varphi d\nu_q,$$ and due to
(\ref{inparts}), it satisfies the property
\begin{equation}\label{inpart}
\int_{U_n}\xi(\psi)\cdot\varphi d\nu_q=\int_{U_n}\psi\cdot
S(\xi)(\varphi) d\nu_q
\end{equation}
for any $\psi\in\mathcal{D}(U_n)'_q$, $ \varphi\in\mathcal{D}(U_n)_q$, and
$\xi\in U_q\mathfrak{su}_{n,1}$. In other words, we have constructed a
morphism of $U_q\mathfrak{su}_{n,1}$-modules\footnote{For an $A$-modules $V$
the dual space $V^*$ is endowed with an $A$-module structure via the
antipode: $\xi f(\cdot)=f(S(\xi)\cdot)$.} $$\mathcal{D}(U_n)'_q\rightarrow
\mathrm{dual}\;\mathrm{of}\;\mathcal{D}(U_n)_q.$$ In fact, this is an
isomorphism. This can be proved by computation just as it was done in
\cite[Section 1]{qd1} in the case $n=1$.

\bigskip

\section{$q$-Berezin transform}\label{3}

The aim of this section is to construct a $q$-analogue of the Berezin
transform in the unit ball. For that purpose we define also $q$-analogues of
the weighted Bergman spaces, Toeplitz operators, and covariant symbols.

\medskip

\subsection{$q$-Weighted Bergman spaces}\label{2.5}

In the classical case the weighted Bergman space is defined as the closure
of the space of holomorphic polynomials with respect to the norm
\begin{equation}\label{weightedmeasure}
\|f\|_\alpha=\left(\frac{\Gamma(n+\alpha+1)}{\Gamma(n+1)\Gamma(\alpha+1)}
\int_{U_n}|f(\mathbf{z})|^2
(1-\|\mathbf{z}\|^2)^{\alpha}dm(\mathbf{z})\right)^{\frac12}
\end{equation}  Here $\alpha>-1$, $dm$ is the normalized Lebesgue
measure, and the multiplier
$\frac{\Gamma(n+\alpha+1)}{\Gamma(n+1)\Gamma(\alpha+1)}$ normalizes the
measure $dm_\alpha=(1-\|\mathbf{z}\|^2)^{\alpha}dm(\mathbf{z})$. Let us
define a $q$-analogue of $dm_\alpha$.

Suppose $\alpha\in\mathbb{Z}_{\geq0}$. The formula (\ref{weightedmeasure}),
along with (\ref{measure}), suggests the following definition of the
$q$-weighted integral
\begin{equation}\label{qwmeasure}
\int_{U_n}fdm_{\alpha,q}{=}\frac{\Gamma_{q^2}(n+\alpha+1)}
{\Gamma_{q^2}(n+1)\Gamma_{q^2}(\alpha+1)} \int_{U_n}f\cdot
y_1^{\alpha-n-1}d\nu_q
\end{equation}
with $y_1{=}1-z_1z_1^*-\ldots -z_nz_n^*\in\mathcal{P}(\mathbb{C}^n)_q$ (see
subsection \ref{2.2}). Here we use the standard notation \cite{GR}:
$$\Gamma_q(x){=}\frac{(q;q)_\infty}{(q^x;q)_\infty} (1-q)^{1-x}$$ with
$$(a;q)_\infty{=}(1-a)(1-aq)(1-aq^2)\ldots,\qquad
(a;q)_\gamma{=}\frac{(a;q)_\infty}{(aq^\gamma;q)_\infty}, \quad
\gamma\in\mathbb{C}.$$ Using (\ref{invint}) and the definition of
$\Gamma_{q^2}$, we may rewrite the integral in the following way
\begin{equation}\label{winvint}
\int_{U_n}fdm_{\alpha,q}=(q^{2\alpha+2};q^2)_n\sum_{\underline{k}\in P(n)
}f_{\underline{0},\underline{0}}(q^{2k_1}, q^{2k_2},\ldots q^{2k_n})\cdot
q^{2k_1(\alpha+1)}\cdot q^{2k_2}\cdot\ldots\cdot q^{2k_n}.
\end{equation}
The latter formula can be used to define the $q$-weighted integral for
arbitrary $\alpha\in\mathbb{R}$.

It is reasonable to define $L^2(dm_{\alpha,q})$ as the completion of
$\mathcal{D}(U_n)_q$ with respect to the norm $$
\|f\|_{\alpha,q}=\left(\int_{U_n} f^*\cdot
fdm_{\alpha,q}\right)^{\frac12}.$$ We shall use the concrete realization of
$L^2(dm_{\alpha,q})$ as the subspace in $\mathcal{D}(U_n)'_q$ of those
distributions for which the right-hand side is finite.

From now on we suppose $\alpha$ is a positive real number. In this case, we
have $\mathcal{P}(\mathbb{C}^n)_q\subset L^2(dm_{\alpha,q})$ (this may be
checked by direct computations using formula (\ref{winvint})). By analogy
with the classical case, we define the $q$-weighted Bergman space
$L_a^2(dm_{\alpha,q})$ as the closure in $L^2(dm_{\alpha,q})$ of the
subspace $\mathbb{C}[\mathbb{C}^n]_q$ of 'holomorphic' polynomials.

\medskip

\begin{proposition}\label{norms}
The monomials $\mathbf{z}^{\underline{m}}$ constitute an orthogonal basis in
$L_a^2(dm_{\alpha,q})$ and $$\|\mathbf{z}^{\underline{m}}\|^2_{\alpha,q}=
\frac{\Gamma_{q^2}(n+\alpha+1)\Gamma_{q^2}(m_1+1)\ldots\Gamma_{q^2}(m_n+1)}
{\Gamma_{q^2}(m_1+\ldots+m_n+n+\alpha+1)}.$$
\end{proposition}

\medskip

\noindent{\bf Proof.} By the above definitions
$$(\mathbf{z}^{\underline{m}},\mathbf{z}^{\underline{l}})_{\alpha,q}
=(q^{2\alpha+2};q^2)_n\sum_{\underline{k}\in P(n) }
f_{\underline{0},\underline{0}}(q^{2k_1}, q^{2k_2},\ldots q^{2k_n})\cdot
q^{2k_1(\alpha+1)}\cdot q^{2k_2}\cdot\ldots\cdot q^{2k_n},$$ where
$f_{\underline{0},\underline{0}}$ is the corresponding term in the expansion
(\ref{form}) of the polynomial
$\mathbf{z}^{*\underline{l}}\mathbf{z}^{\underline{m}}$. Obviously, the
$f_{\underline{0},\underline{0}}$-term vanishes for $\underline{m}\ne
\underline{l}$. This implies the pairwise orthogonality of the monomials.
Suppose $\underline{m}=\underline{l}$. It is easy to show that
\begin{equation}\label{sum}
\mathbf{z}^{*\underline{m}}\mathbf{z}^{\underline{m}}=y_2^{m_1}y_3^{m_2}\ldots
y_n^{m_{n-1}}(q^2\frac{y_1}{y_2};q^2)_{m_1}(q^2\frac{y_2}{y_3};q^2)_{m_2}
\ldots (q^2\frac{y_{n-1}}{y_n};q^2)_{m_{n-1}} (q^2y_n;q^2)_{m_n}.
\end{equation}
Denote by $f(y_1,y_2,\ldots y_n)$ the right-hand side of the latter
equality. One has {\small$$\sum_{\underline{k}\in P(n)}f(q^{2k_1},
q^{2k_2},\ldots q^{2k_n})\cdot q^{2k_1(\alpha+1)+2k_2+\ldots+2k_n}=
\frac1{(1-q^2)^n}\prod_{j=0}^{n-1}\int_0^1t^{\alpha+m_1+\ldots+m_j+j}
(q^2t;q^2)_{m_{j+1}}d_{q^2}t$$} where $\int_0^1f(t)d_{q^2}t$ is the Jackson
integral given by
\begin{equation}\label{jackson}\int_0^1f(t)d_{q^2}t=(1-q^2)\sum_{l=0}^\infty
f(q^{2l})q^{2l}.\end{equation} What remains is to use the following well
known formula \cite{GR}: $$\int_0^1t^{a}(q^2t;q^2)_{b}d_{q^2}t=
\frac{\Gamma_{q^2}(a+1)\Gamma_{q^2}(b+1)}{\Gamma_{q^2}(a+b+2)}.$$ \hfill
$\blacksquare$

\medskip

In the classical case the Hilbert space $L^2(dm_{\alpha})$ admits a natural
unitary $\widetilde{SU(n,1)}$-action, where $\widetilde{SU(n,1)}$ is the
universal covering of $SU(n,1)$. The invariant subspace $L_a^2(dm_{\alpha})$
of holomorphic functions is called the representation of the holomorphic
discrete series \cite{Wallach}, \cite{Arazy}. The corresponding
infinitesimal $U\mathfrak{su}_{n,1}$-action is obtained by a simple twisting
of the natural action. Below, we produce an analogue of the twisted
$U\mathfrak{su}_{n,1}$-action.

It was proved in \cite[Section 6]{SSV2} that there exists a unique
representation $\pi_\alpha$ of $U_q\mathfrak{su}_{n,1}$ in
$\mathcal{P}(\mathbb{C}^n)_q$ and $\mathcal{D}(U_n)_q$ such that for all $f
\in\mathcal{P}(\mathbb{C}^n)_q$ (or $f\in\mathcal{D}(U_n)_q$)
\begin{equation}\label{pl1}\pi_\alpha(E_j)f=
\begin{cases} E_jf, & j \ne n, \\ E_nf-q^{1/2}\dfrac{1-q^{2
\alpha+2n+2}}{1-q^2}(K_nf)z_n, & j=n,
\end{cases}
\end{equation}
\begin{equation}\pi_\alpha(F_j)f=
\begin{cases} F_jf, & j \ne n, \\ q^{-\alpha-n-1} F_nf, & j=n,
\end{cases}
\end{equation}
\begin{equation}\label{pl3}\pi_\alpha(K_j^{\pm 1})f=
\begin{cases} K_j^{\pm 1}f, & j \ne n, \\ q^{\pm(\alpha+n+1)}K_n^{\pm 1}f,
& j=n,
\end{cases}
\end{equation}
where $E_jf, F_jf, K_j f$ are defined in Proposition \ref{2.1} ((\ref{ef})
and (\ref{hfe}) for $f\in\mathcal{D}(U_n)_q$). In order to not consider the
polynomials and finite functions separately, it will be convenient sometimes
to regard $\pi_\alpha$ as a representation of $U_q\mathfrak{su}_{n,1}$ in
$\mathcal{F}(U_n)_q$.

The following statement is proved in \cite[Section 6]{SSV2}.

\medskip

\begin{proposition}\label{unitary}
For any $\varphi_1,\varphi_2\in\mathcal{F}(U_n)_q$ and any $\xi\in
U_q\mathfrak{su}_{n,1}$ $$ (\pi_\alpha(\xi)\varphi_1,\varphi_2)_{\alpha,q}=
(\varphi_1,\pi_\alpha(\xi^*)\varphi_2)_{\alpha,q}$$ with
$(\cdot\;,\cdot\;)_{\alpha,q}$ being the scalar product in
$L^2(dm_{\alpha,q})$.
\end{proposition}

\medskip

Note that the subspace $\mathbb{C}[\mathbb{C}^n]_q\subset
\mathcal{F}(U_n)_q$ is actually a $U_q\mathfrak{su}_{n,1}$-submodule. It
should be treated as a module of the holomorphic discrete series for
$U_q\mathfrak{su}_{n,1}$.

Before we proceed further, let us fix more notation. The space
$\mathbb{C}[\mathbb{C}^n]_q$ (as well as $\mathcal{F}(U_n)_q$) admits two
different $U_q\mathfrak{su}_{n,1}$-actions, namely, the original one which
was introduced earlier in Proposition \ref{2.1} and the twisted one which
appeared in this subsection. In order to emphasize that a space is
considered together with the twisted $U_q\mathfrak{su}_{n,1}$-action
$\pi_\alpha$ we shall add the subscript $\alpha$ in the notation, for
example, $\mathbb{C}[\mathbb{C}^n]_{q,\alpha}$. The subscript will also
indicate the pre-Hilbert space structure given by the scalar product
$(\cdot,\cdot)_{\alpha,q}$.

The following important observation relates the usual and the twisted
$U_q\mathfrak{su}_{n,1}$-actions: the multiplication map
$\mathcal{F}(U_n)_q\otimes\mathcal{F}(U_n)_{q,\alpha}\rightarrow\mathcal{F}
(U_n)_{q,\alpha}$ is a morphism of $U_q\mathfrak{su}_{n,1}$-modules. This
can be proved by direct computations.

\medskip

\subsection{Toeplitz operators}

In this subsection we produce $q$-analogues of Toeplitz operators with
polynomial and finite symbols. As in the previous subsection, it is
convenient to consider symbols from the algebra $\mathcal{F}(U_n)_q$ when
there is no necessity to consider the cases of polynomial and finite symbols
separately.

Let $P_{\alpha,q}$ be the orthogonal projection in $L^2(dm_{\alpha,q})$ onto
$L_a^2(dm_{\alpha,q})$. The Toeplitz operator $T_f$ with the symbol
$f\in\mathcal{F}(U_n)_q$ is defined as follows
$$T_f:\mathbb{C}[\mathbb{C}^n]_{q,\alpha}\rightarrow\mathbb{C}
[\mathbb{C}^n]_{q,\alpha},\quad T_f(\psi)=P_{\alpha,q}(f\cdot\psi).$$

To formulate the principal result of the present subsection, we recall the
following well known construction. Let $A$ be a Hopf algebra and $V$ an
$A$-module. Then the space ${\rm End}(V)$ admits the following 'canonical'
structure of $A$-module: for $\xi\in A$, $T\in{\rm End}(V)$ $$
\xi(T)=\sum_j\xi'_j\cdot T\cdot S(\xi''_j),$$ where $S$ is the antipode of
$A$, $\Delta(\xi)=\sum_j\xi'_j\otimes\xi''_j$ (with $\Delta$ being the
comultiplication), and the elements in the right-hand side are multiplied
within the algebra ${\rm End}(V)$. This action of $A$ in ${\rm End}(V)$
makes ${\rm End}(V)$ into an $A$-module algebra.

\medskip

\begin{proposition}\label{toepinv} The linear map
$$\mathcal{F}(U_n)_q\rightarrow\mathrm{End}(\mathbb{C}[\mathbb{C}^n]
_{q,\alpha}),\quad f\mapsto T_f$$ is a morphism of
$U_q\mathfrak{su}_{n,1}$-modules.
\end{proposition}

\medskip

\noindent{\bf Proof.} In the classical case the projection $P_{\alpha,q}$
intertwines the twisted $\widetilde{SU(n,1)}$-action in $L^2(dm_{\alpha})$
and $L_a^2(dm_{\alpha})$. In the quantum case the
$U_q\mathfrak{su}_{n,1}$-action $\pi_\alpha$ is not defined on the entire
$L^2(dm_{\alpha,q})$. Nevertheless, the intertwining property may be
formulated due to the equality
$$P_{\alpha,q}(\mathcal{F}(U_n)_{q,\alpha})=
\mathbb{C}[\mathbb{C}^n]_{q,\alpha}.$$ To prove the equality, we endow the
space $\mathbb{C}[\mathbb{C}^n]_{q,\alpha}$ with the obvious
$\mathbb{Z}_{\geq0}$-grading by powers of monomials and observe that each
$f\in\mathcal{F}(U_n)_{q,\alpha}$ is orthogonal to all but finitely many
homogeneous components of $\mathbb{C}[\mathbb{C}^n]_{q,\alpha}$. Now it is
evident that the operator
$P_{\alpha,q}:\mathcal{F}(U_n)_{q,\alpha}\rightarrow\mathbb{C}
[\mathbb{C}^n]_{q,\alpha}$ is a morphism of
$U_q\mathfrak{su}_{n,1}$-modules: this follows from
$U_q\mathfrak{su}_{n,1}$-invariance of the orthogonal complement of
$\mathbb{C}[\mathbb{C}^n]_{q,\alpha}$ in $\mathcal{F}(U_n)_{q,\alpha}$.

\medskip

\noindent{\bf Remark.} One can use the above observation to show that the
projection $P_{\alpha,q}$ is a $q$-integral operator with a simple kernel
which is a $q$-analogue of the so-called Bergman kernel; see \cite{SSV2}.
But this result is not needed in the present paper.

\medskip

Let us get back to proving the proposition. Let $f\in\mathcal{F}(U_n)_q$.
Denote by $\hat f= \hat{f}_\alpha$ the endomorphism of
$\mathcal{F}(U_n)_{q,\alpha}$ given by
$\hat{f}_\alpha(\varphi)=f\cdot\varphi$. The map
$$\mathcal{F}(U_n)_q\rightarrow\mathrm{End}(\mathcal{F}(U_n)_{q,\alpha}),
\quad f\mapsto \hat{f}_\alpha$$ is a morphism
$U_q\mathfrak{su}_{n,1}$-modules (see the remark at the end of the previous
subsection). It remains to use the equality
$T_f=P_{q,\alpha}\cdot\hat{f}_\alpha\mid_{\mathbb{C}
[\mathbb{C}^n]_{q,\alpha}}$. \hfill$\blacksquare$

\medskip

\subsection{On the case of polynomial symbols}

For purposes of this paper we need mostly Toeplitz operators with finite
symbols. However, we present some results about the polynomial case which
seem to be interesting by themselves.

Let $\hat{z}_i$ stand for the operator in
$\mathbb{C}[\mathbb{C}^n]_{q,\alpha}$ of the left multiplication by $z_i$,
$\hat{z}^*_i$ for the adjoint operator. The following observation is
straightforward: for a polynomial $f=\sum
a_{\underline{l},\underline{m}}\mathbf{z}^{*\underline{l}}
\mathbf{z}^{\underline{m}}$
one has $T_f=\sum
a_{\underline{l},\underline{m}}\hat{\mathbf{z}}^{*\underline{l}}
\hat{\mathbf{z}}^{\underline{m}}$. Thus all Toeplitz operators with
polynomial symbols belong to the unital subalgebra in
$\mathrm{End}(\mathbb{C}[\mathbb{C}^n]_{q,\alpha})$ generated by $\hat{z}_i$
and $\hat{z}^*_i$, $i=1,2,\ldots n$. Denote this subalgebra by
$\mathcal{P}_{n,q,\alpha}$. It is an involutive algebra with the involution
given by $*:\hat{z}_i\mapsto\hat{z}^*_i$.

\medskip

\begin{proposition}
The $U_q\mathfrak{su}_{n,1}$-action in
$\mathrm{End}(\mathbb{C}[\mathbb{C}^n]_{q,\alpha})$ induces a
$U_q\mathfrak{su}_{n,1}$-module algebra structure in
$\mathcal{P}_{n,q,\alpha}$.
\end{proposition}

\medskip

\noindent{\bf Proof.} First we have to establish
$U_q\mathfrak{su}_{n,1}$-invariance of the subspace
$\mathcal{P}_{n,q,\alpha}$ in
$\mathrm{End}(\mathbb{C}[\mathbb{C}^n]_{q,\alpha})$, that is, to show that
$\xi(\hat{z}^{*a}_i\hat{z}^b_j)$ belong to $\mathcal{P}_{n,q,\alpha}$ for
any $a,b\in\mathbb{Z}_{\geq0}$, $i,j=1,2,\ldots n$, and $\xi\in
U_q\mathfrak{su}_{n,1}$. According to proposition \ref{toepinv}
$$\xi(\hat{z}^{*a}_i\hat{z}^b_j)=\xi(
T_{z^{*a}_iz^b_j})=T_{\xi(z^{*a}_iz^b_j)}.$$ By the remark, preceding the
present proposition, $T_{\xi(z^{*a}_iz^b_j)}\in\mathcal{P}_{n,q,\alpha}$. It
remains to prove module algebra property (\ref{agree}) for
$\mathcal{P}_{n,q,\alpha}$. But it may be derived easily from Proposition
\ref{unitary}.\hfill$\blacksquare$

\medskip

\begin{proposition}\label{gkl} The following commutation relations hold
$$ \hat{z}_i\hat{z}_j=q\hat{z}_j\hat{z}_i, \quad i<j, $$
$$ \hat{z}_i^*
\hat{z}_j=q\hat{z}_j\hat{z}_i^*+q^{2\alpha+2n}\left(\hat{z}_i^*
\hat{z}_j(q^2+(1-q^2)\sum_{k=1}^n\hat{z}_k\hat{z}_k^*)-q\hat{z}_j
\hat{z}_i^*\right),
\quad i\ne j,$$
\begin{equation*}
\begin{split}
 \hat{z}_j^*\hat{z}_j&=q^2\hat{z}_j\hat{z}_j^*+(1-q^{2})
(1-\sum_{k=j+1}^n\hat{z}_k\hat{z}_k^*)\\ &\qquad +q^{2\alpha+2n}(\hat{z}_j^*
\hat{z}_j(q^2+(1-q^2)\sum_{k=1}^n\hat{z}_k\hat{z}_k^*)-q^2\hat{z}_j
\hat{z}_j^*-(1-q^{2})
\sum_{k=1}^j\hat{z}_k\hat{z}_k^*).
\end{split}
\end{equation*}
\end{proposition}

\medskip

\noindent{\bf Proof.} The relations may be deduced from the formulas
$$\hat{z}_i(\mathbf{z}^{\underline{m}})=q^{m_n+m_{n-1}+\ldots+m_{i+1}}
\mathbf{z}^{\underline{m'}}$$ (where $\underline{m'}=(m_1,\ldots,
m_i+1,\ldots,m_n)$),

$$\hat{z}^*_i(\mathbf{z}^{\underline{m}})=0, \quad m_i=0, $$

$$\hat{z}^*_i(\mathbf{z}^{\underline{m}})=\frac{q^{m_n+m_{n-1}+\ldots+
m_{i+1}}(1-q^{2m_i})}
{1-q^{2m_1+2m_2+\ldots+2m_n+2n+2\alpha}}\mathbf{z}^{\underline{m}''}$$
(where $\underline{m''}=(m_1,\ldots, m_i-1,\ldots,m_n)$). The first equality
is obvious, the other follow from Proposition \ref{norms}.
\hfill$\blacksquare$

\medskip
Note that for $\alpha=\infty$ the above relations coincides with the defining
relations for $\mathcal{P}(\mathbb{C}^n)_q$.

The unital involutive algebra given by the above commutation relations can
be viewed as a two-parameter deformation of the polynomial algebra on
$\mathbb{C}^n$. It is easy to show that for $n=1$ the algebra is isomorphic
to the one considered in details in \cite{KL} (see the Introduction).

\medskip

\subsection{Covariant symbols and $q$-Berezin transform}\label{3.3}

In this subsection, we define the notion of covariant symbols of operators
on the $q$-weighted Bergman spaces and use it to define a $q$-analogue of
the Berezin transform.

To define covariant symbols we need certain inner product in the space of
operators on a $q$-weighted Bergman space (see the Introduction). This
product is defined via the so-called $q$-trace. So we recall first its
definition.

Let $V$ be a finite-dimensional $U_q\mathfrak{sl}_{n+1}$-module. The
$q$-trace is the linear functional on $\mathrm{End}(V)$ given by
\begin{equation}\label{qtrace}
\mathrm{tr}_q:T\mapsto
\mathrm{tr}\left(T\cdot\prod_{j=1}^nK_j^{-j(n+1-j)}\right).\end{equation}
The following well known observation explains the importance of this
functional (see \cite{ChP}): the $q$-trace is an invariant linear functional
on $\mathrm{End}(V)$, i.e.
$\mathrm{tr}_q(\xi(T))=\varepsilon(\xi)\cdot\mathrm{tr}_q(T)$ for any
$\xi\in U_q\mathfrak{sl}_{n+1}$ and $T\in\mathrm{End}(V)$.

 We will modify the definition slightly
for the infinite-dimensional space $V=\mathbb{C}[\mathbb{C}^n]_{q,\alpha}$.
Let $\mathrm{End}_{0}(\mathbb{C}[\mathbb{C}^n]_{q,\alpha})$ be the subspace
in $\mathrm{End}(\mathbb{C}[\mathbb{C}^n]_{q,\alpha})$ of those
automorphisms whose matrices in the basis $\{\mathbf{z}^{\underline{m}}\}$
have only finitely many non-zero entries. The $q$-trace (\ref{qtrace}) is a
well defined linear functional on
$\mathrm{End}_{0}(\mathbb{C}[\mathbb{C}^n]_{q,\alpha})$. The invariance
property still holds in this case.

\medskip

\begin{proposition}\label{qtr1}
The $q$-trace is an invariant linear functional on
$\mathrm{End}_{0}(\mathbb{C}[\mathbb{C}^n]_{q,\alpha})$.
\end{proposition}

\medskip

\noindent{\bf Sketch of a proof.} A standard proof in the finite-dimensional
case uses the canonical isomorphism of $U_q\mathfrak{sl}_{n+1}$-modules
$\mathrm{End}(V)\simeq V\otimes V^*$ with $V^*$ being the dual
$U_q\mathfrak{sl}_{n+1}$-module \cite{ChP}. The statement of Proposition
\ref{qtr1} can be proved by the same argument since
$\mathrm{End}_{0}(\mathbb{C}[\mathbb{C}^n]_{q,\alpha})\simeq
\mathbb{C}
[\mathbb{C}^n]_{q,\alpha}\otimes\mathbb{C}[\mathbb{C}^n]_{q,\alpha}^*$.
\hfill$\blacksquare$

\medskip

\noindent{\bf Remarks.} i) The formula (\ref{qtrace}) is the same as the one
which defines the invariant integral on the quantum ball (see Proposition
\ref{integral}). But Proposition \ref{integral} does not follow formally
from the above statement since the $\mathcal{D}(U_n)_q$-module $H$ is not a
$U_q\mathfrak{sl}_{n+1}$-module.

ii) Using the equality (\ref{inpart2}) below and the same idea as in the
proof of Proposition \ref{scpro} one can show that the scalar product in
$\mathrm{End}_{0}(\mathbb{C}[\mathbb{C}^n]_{q,\alpha})$ given by
$(T_1,T_2)=\mathrm{tr}_q(T_2^*\cdot T_1)$ is
$U_q\mathfrak{su}_{n,1}$-invariant. The unitary representation of
$U_q\mathfrak{su}_{n,1}$ in
$\mathrm{End}_{0}(\mathbb{C}[\mathbb{C}^n]_{q,\alpha})$ should be treated as
a $q$-analogue of the canonical representation of $SU(n,1)$ \cite{Dijk}.

\medskip

For a particular $\alpha$, we shall use the modified $q$-trace
$\mathrm{Tr}_q$ on $\mathrm{End}_{0}(\mathbb{C}[\mathbb{C}^n]_{q,\alpha})$
which differs from $\mathrm{tr}_q$ by a constant:
$$\mathrm{Tr}_q=
q^{n(n+\alpha+1)}\frac{(q^2;q^2)_n}{(q^{2+2\alpha};q^2)_n}\mathrm{tr}_q.$$
It is easy to deduce the following explicit formula for $\mathrm{Tr}_q$. Let
$A\in\mathrm{End}_{0}(\mathbb{C}[\mathbb{C}^n]_{q,\alpha})$ be the
endomorphism given by $A(\mathbf{z}^{\underline{m}})=
\sum_{\underline{k}}A_{\underline{k}}^{\underline{m}}\cdot
\mathbf{z}^{\underline{k}}$. Then
\begin{equation}\label{qtra}
\mathrm{Tr}_q(A)=
\frac{(q^2;q^2)_n}{(q^{2+2\alpha};q^2)_n}\cdot\sum_{\underline{m}}
A_{\underline{m}}^{\underline{m}}\cdot q^{-2n(m_1+m_2+\ldots+m_n)}\cdot
q^{2m_2+\ldots+2m_n}\cdot\ldots\cdot q^{2m_n}.
\end{equation}

Now we are ready to define the notion of covariant symbols. Let
$T\in\mathrm{End}(\mathbb{C}[\mathbb{C}^n]_{q,\alpha})$. A distribution
$\sigma(T)\in\mathcal{D}(U_n)'_q$ is said to be the covariant symbol of $T$
if for any $f\in\mathcal{D}(U_n)_q$
\begin{equation}\label{cov}
\int_{U_n}\sigma(T)\cdot fd\nu_q=\mathrm{Tr}_q(T\cdot T_f).
\end{equation}
Note that for $f\in\mathcal{D}(U_n)_q$ one has
$T_f\in\mathrm{End}_0(\mathbb{C}[\mathbb{C}^n]_{q,\alpha})$. Thus the
right-hand side of (\ref{cov}) is well defined. The arguments cited at the
end of subsection 2.4 imply existence and uniqueness of the covariant symbol
of an arbitrary endomorphism $F$. Actually, the map $T\mapsto\sigma(T)$ is
conjugated to the map $f\mapsto T_f$.

We have the following elementary property of the $q$-trace, which can be
proved by the same arguments as the formula of integrating by parts (see the
proof of Proposition \ref{scpro}): for any
$T\in\mathrm{End}(\mathbb{C}[\mathbb{C}^n]_{q,\alpha})$,
$T_0\in\mathrm{End}_0(\mathbb{C}[\mathbb{C}^n]_{q,\alpha})$, and $\xi\in
U_q\mathfrak{su}_{n,1}$

\begin{equation}\label{inpart2}
\mathrm{Tr}_q(\xi(T)\cdot T_0)=\mathrm{Tr}_q(T\cdot S(\xi)(T_0)).
\end{equation}
The formula implies

\medskip

\begin{proposition}\label{invcov}
The map $$\mathrm{End}(\mathbb{C}[\mathbb{C}^n]_{q,\alpha})\rightarrow
\mathcal{D}(U_n)'_q,\qquad T\mapsto\sigma(T)$$ is a morphism of
$U_q\mathfrak{su}_{n,1}$-modules.
\end{proposition}

\medskip

\noindent{\bf Proof.} The equality (\ref{inpart2}) gives an identification
of the $U_q\mathfrak{su}_{n,1}$-module
$\mathrm{End}(\mathbb{C}[\mathbb{C}^n]_{q,\alpha})$ with the dual of
$\mathrm{End}_0(\mathbb{C}[\mathbb{C}^n]_{q,\alpha})$ (see the definition of
the dual module in subsection \ref{2.4}). What remains is to use Proposition
\ref{toepinv} and the observation that the map
$\mathrm{End}(\mathbb{C}[\mathbb{C}^n]_{q,\alpha})\rightarrow
\mathcal{D}(U_n)'_q$, $T\mapsto\sigma(T)$ is conjugated to
$\mathcal{D}(U_n)_q\rightarrow
\mathrm{End}_0(\mathbb{C}[\mathbb{C}^n]_{q,\alpha})$,
$ f\mapsto T_f$. \hfill$\blacksquare$

\medskip

We are in position to define the $q$-Berezin transform $B_{q,\alpha}$. It is
defined as the linear map from $\mathcal{D}(U_n)_q$ to $\mathcal{D}(U_n)'_q$
which sends a finite function to the covariant symbol of the corresponding
Toeplitz operator:
$$B_{q,\alpha}:
f\mapsto\sigma(T_f).$$ The following crucial statement is
straightforward.

\medskip

\begin{proposition}\label{berinv}
The $q$-Berezin transform is a morphism of $U_q\mathfrak{su}_{n,1}$-modules.
\end{proposition}

\medskip

Due to Proposition \ref{f00}, any morphism of
$U_q\mathfrak{su}_{n,1}$-modules
$T:\mathcal{D}(U_n)_q\rightarrow\mathcal{D}(U_n)'_q$ is completely
determined by the element $T(f_0)\in\mathcal{D}(U_n)'_q$. Thus it would be
very useful to compute $B_{q,\alpha}(f_0)$.

\medskip

\begin{proposition}\label{bf0}
$$B_{q,\alpha}(f_0)=(q^{2\alpha+2};q^2)_n\cdot y_1^{\alpha+n+1}.$$
\end{proposition}

\medskip

\noindent{\bf Proof.} We have to check that for any $f\in\mathcal{D}(U_n)_q$
\begin{equation}\label{equ}
(q^{2\alpha+2};q^2)_n\int_{U_n}y_1^{\alpha+n+1}\cdot
fd\nu_q=\mathrm{Tr}_q(T_{f_0}\cdot T_f).
\end{equation}
Let us denote by $f^{\underline{k}}_{\underline{i},\underline{j}}$
($\underline{k}\in P(n)$, $\underline{i},
\underline{j}\in\mathbb{Z}^n_{\geq0}$,
$\underline{i}\times\underline{j}=\underline{0}$) the finite function given
by
$f^{\underline{k}}_{\underline{i},\underline{j}}=\mathbf{z}^{\underline{i}}
f_{\underline{k}}(y_1,y_2,\ldots y_n)\mathbf{z}^{*\underline{j}}$ with $
f_{\underline{k}}(q^{2l_1}, q^{2l_2},\ldots q^{2l_n})=\left
\{\begin{array}{ccl}1, & \underline{k}=\underline{l}
\\ 0, & {\rm otherwise.}\end{array}\right.$
For example, $f^{\underline{0}}_{\underline{0},\underline{0}}=f_0$.
Evidently, the functions $f^{\underline{k}}_{\underline{i},\underline{j}}$
constitute a basis in $\mathcal{D}(U_n)_q$.

It follows from the definitions that the both hand sides of (\ref{equ})
vanish for $f:=f^{\underline{k}}_{\underline{i},\underline{j}}$ provided
$\underline{i}\ne\underline{0}$ or $\underline{j}\ne\underline{0}.$ Thus we
have to verify (\ref{equ}) for $f:=f_{\underline{k}}$. It follows from
(\ref{invint}) that $$ (q^{2\alpha+2};q^2)_n\int_{U_n}y_1^{\alpha+n+1}\cdot
f_{\underline{k}}d\nu_q=(q^{2\alpha+2};q^2)_n\cdot(q^2;q^2)_n\cdot
q^{2k_1(\alpha+1)}\cdot q^{2k_2}\cdot\ldots\cdot q^{2k_n}. $$ Let us compute
the right-hand side of (\ref{equ}) with $f:=f_{\underline{k}}$. Evidently,
$$f_{0}\cdot \mathbf{z}^{\underline{m}}=\left \{\begin{array}{ccl}f_{0}, &
\underline{m}=\underline{0}
\\ 0, & {\rm otherwise.}\end{array}\right.$$
Besides, $\int_{U_n}f_{0}dm_{\alpha,q}=(q^{2\alpha+2};q^2)_n.$ These
equalities mean that $$T_{f_{0}} (\mathbf{z}^{\underline{m}})=\left
\{\begin{array}{ccl} (q^{2\alpha+2};q^2)_n,& \underline{m}=\underline{0}
\\ 0,& {\rm otherwise.}\end{array}\right.$$
Clearly, the vector $f_{\underline{k}}\cdot \mathbf{z}^{\underline{m}}\in
L^2(dm_{\alpha,q})$ is orthogonal to any monomial
$\mathbf{z}^{\underline{l}}$ except the case $\underline{l}=\underline{m}$.
In particular, the monomials are eigenvectors of $T_{f_{\underline{k}}}$.
Thus $$\mathrm{Tr}_q(T_{f_0}\cdot T_{f_{\underline{k}}})=(q^2;q^2)_n
(T_{f_{\underline{k}}}(1),1)_{\alpha,q}=(q^2;q^2)_n\int_{U_n}f
_{\underline{k}} dm_{\alpha,q}=$$ $$
=(q^{2\alpha+2};q^2)_n\cdot(q^2;q^2)_n\cdot q^{2k_1(\alpha+1)}\cdot
q^{2k_2}\cdot\ldots\cdot q^{2k_n}.$$ \hfill$\blacksquare$

\bigskip

\section
{$q$-Laplace-Beltrami operator and
associated $q$-spherical transform}\label{4}

In this section we define a $q$-analogue of the $SU(n.1)$-invariant Laplace
operator (the Laplace-Beltrami operator) on the unit ball and study a
$q$-analogue of the spherical transform in the unit ball. Namely, we
calculate the 'radial part' of the $q$-Laplace-Beltrami operator, find
$q$-spherical functions, and present an inversion formula for the
$q$-spherical transform.

\medskip

\subsection{$q$-Laplace-Beltrami operator}\label{4.1}

In this subsection, we consider the asymptotic expansion of the $q$-Berezin
transform $B_{q,\alpha}$ at the limit $t=q^{2\alpha}\to0$. A $q$-analogue of
the Laplace-Beltrami operator on the ball is defined as the coefficient at
$t$ in that expansion.

Recall (see the proof of Proposition \ref{bf0}) the notation
$f^{\underline{k}}_{\underline{i},\underline{j}}$ ($\underline{k}\in P(n)$,
$\underline{i}, \underline{j}\in\mathbb{Z}^n_{\geq0}$,
$\underline{i}\times\underline{j}=\underline{0}$). Suppose $T$ is the linear
map from $\mathcal{D}(U_n)_q$ to $\mathcal{D}(U_n)'_q$. The numbers

$$T_{\underline{i},\underline{j};\underline{p},\underline{r}}^
{\underline{k},\underline{s}}=\int_{U_n}T
(f^{\underline{k}}_{\underline{i},\underline{j}})\cdot
f^{\underline{s}}_{\underline{p},\underline{r}}d\nu_q.$$ will be called the
matrix entries of $T$. Due to non-degeneracy of the pairing
$\mathcal{D}(U_n)'_q\times\mathcal{D}(U_n)_q\rightarrow\mathbb{C}$
(subsection \ref{2.4}), the matrix entries determine $T$ completely.

Let us denote the matrix entries of the $q$-Berezin transform $B_{q,\alpha}$
by $B_{\underline{i},\underline{j};\underline{p},\underline{r}}^
{\underline{k},\underline{s}}(\alpha)$. Introduce the new variable
$t=q^{2\alpha}$. We will regard the matrix entries
$B_{\underline{i},\underline{j};\underline{p},\underline{r}}^
{\underline{k},\underline{s}}(\alpha)$ as functions of $t$ and use the
notation $B_{\underline{i},\underline{j};\underline{p},\underline{r}}^
{\underline{k},\underline{s}}(t):=
B_{\underline{i},\underline{j};\underline{p},
\underline{r}}^ {\underline{k},\underline{s}}(\alpha)$,
$B_{q,t}:=B_{q,\alpha}$.

\medskip

\begin{proposition}\label{series}
There exists a sequence $\{B_{q,j}\}_{j\in\mathbb{N}}$ of linear
endomorphisms of $\mathcal{D}(U_n)_q$ such that

i) $B_{q,j}$ are independent of $t$ (i.e. their matrix entries are
independent of $t$);

ii) each $B_{q,j}$ is a $U_q\mathfrak{su}_{n,1}$-module morphism;

iii) for any $f\in\mathcal{D}(U_n)_q$
\begin{equation}\label{ex}
B_{q,t}(f)=f+\sum_{j=1}^\infty B_{q,j}(f)\cdot t^j\end{equation}
where the series is convergent in $\mathcal{D}(U_n)'_q$.
\end{proposition}

\medskip

\noindent{\bf Proof.} First of all, we will construct a sequence
$\{B_{q,j}\}_{j\in\mathbb{Z}_{\geq0}}$ of linear operators from
$\mathcal{D}(U_n)_q$ to $\mathcal{D}(U_n)'_q$ such that each $B_{q,j}$ is
independent of $t$ and $ B_{q,t}(f)=\sum_{j=0}^\infty B_{q,j}(f)\cdot t^j$
for any $f\in\mathcal{D}(U_n)_q$ (the series is convergent in the sense that
all the corresponding series of matrix entries are convergent). For that
purpose we need

\medskip

\begin{lemma}
The matrix entries
$B_{\underline{i},\underline{j};\underline{p},\underline{r}}^
{\underline{k},\underline{s}}(t)$ are polynomials in $t$.
\end{lemma}

\medskip

\noindent{\bf Proof of the lemma.} Due to Proposition \ref{f00},
for any
 $\underline{k},
\underline{i}, \underline{j}\in\mathbb{Z}^n_{\geq0}$
($\underline{i}\times\underline{j}=\underline{0}$) there exists
$\xi^{\underline{k}}_{\underline{i},\underline{j}}\in
U_q\mathfrak{su}_{n,1}$ such that
$f^{\underline{k}}_{\underline{i},\underline{j}}=
\xi^{\underline{k}}_{\underline{i},\underline{j}}(f_0).$ Then, $$
B_{\underline{i},\underline{j};\underline{p},\underline{r}}^
{\underline{k},\underline{s}}(t)=\int_{U_n}B_{q,t}
(f^{\underline{k}}_{\underline{i},\underline{j}})\cdot
f^{\underline{s}}_{\underline{p},\underline{r}}d\nu_q=\int_{U_n}B_{q,t}
(\xi^{\underline{k}}_{\underline{i},\underline{j}}(f_0))\cdot
f^{\underline{s}}_{\underline{p},\underline{r}}d\nu_q=\int_{U_n}B_{q,t}
(f_0)\cdot S(\xi^{\underline{k}}_{\underline{i},\underline{j}})
(f^{\underline{s}}_{\underline{p},\underline{r}})d\nu_q,$$ where the last
equality is due to Proposition \ref{berinv} and the equality (\ref{inpart}).
Thus to prove the lemma it suffices to establish that the entries $
B_{\underline{0},\underline{0};\underline{p},\underline{r}}^
{\underline{0},\underline{s}}(t)=\int_{U_n}B_{q,t} (f_0)\cdot
f^{\underline{s}}_{\underline{p},\underline{r}}d\nu_q$ are polynomials.

Recall (Proposition \ref{bf0}) that $B_{q,\alpha}(f_0)=
(q^{2\alpha+2};q^2)_n\cdot y_1^{\alpha+n+1}.$ Thus $$\int_{U_n}B_{q,t}
(f_0)\cdot
f^{\underline{s}}_{\underline{p},\underline{r}}d\nu_q=(q^{2\alpha+2};q^2)_n
\int_{U_n}y_1^{\alpha+n+1}\cdot
f^{\underline{s}}_{\underline{p},\underline{r}}d\nu_q.$$ The latter integral
vanishes with $\underline{p}\ne\underline{0}$ or
$\underline{r}\ne\underline{0}.$ Consider the integral
$(q^{2\alpha+2};q^2)_n \int_{U_n}y_1^{\alpha+n+1}\cdot
f_{\underline{k}}d\nu_q$ ($f_{\underline{k}}$ are defined in the proof of
Proposition \ref{bf0}). Obviously,
\begin{equation*}
\begin{split}
(q^{2\alpha+2};q^2)_n \int_{U_n}y_1^{\alpha+n+1}\cdot
f_{\underline{k}}d\nu_q&=
(q^{2\alpha+2};q^2)_n(q^{2};q^2)_nq^{2k_1(\alpha+1)}q^{2k_2}\ldots
q^{2k_n}=\\ &=(q^2t;q^2)_n(q^{2};q^2)_nq^{2k_1}t^{k_1} q^{2k_2}\ldots
q^{2k_n}.
\end{split}
\end{equation*}
\hfill$\blacksquare$

\medskip

Now we can define $\{B_{q,m}\}_{m\in\mathbb{Z}_{\geq0}}$ by their matrix
coefficients
$$(B_{q,m})_{\underline{i},\underline{j};\underline{p},\underline{r}}^
{\underline{k},\underline{s}}=\frac1{m!}\frac{d^mB_{\underline{i},
\underline{j};\underline{p},\underline{r}}^
{\underline{k},\underline{s}}(t)}{dt^m}\mid_{t=0}.$$ This definition implies
the $U_q\mathfrak{su}_{n,1}$-invariance of the maps $B_{q,m}$. Indeed, we
should prove coincidence of the maps $B_{q,m}\cdot\xi$ and $\xi\cdot
B_{q,m}$ for any $\xi\in U_q\mathfrak{su}_{n,1}$, or, equivalently, equality
of all matrix entries
$(B_{q,m}\cdot\xi)_{\underline{i},\underline{j};\underline{p},\underline{r}}^
{\underline{k},\underline{s}}$ and $(\xi\cdot
B_{q,m})_{\underline{i},\underline{j};\underline{p},\underline{r}}^
{\underline{k},\underline{s}}$. Due to $U_q\mathfrak{su}_{n,1}$-invariance
of $B_{q,t}$, one has
$$(B_{q,t}\cdot\xi)_{\underline{i},\underline{j};\underline{p},
\underline{r}}^ {\underline{k},\underline{s}}=(\xi\cdot
B_{q,t})_{\underline{i},\underline{j};\underline{p},\underline{r}}^
{\underline{k},\underline{s}}$$ for any
$\underline{i},\underline{j},\underline{p},\underline{r},
\underline{k},\underline{s}.$ What remains is to differentiate the latter
equality.

Let us show that the image of the map $B_{q,m}:
\mathcal{D}(U_n)_q\rightarrow\mathcal{D}(U_n)'_q$ is contained in
$\mathcal{D}(U_n)_q$. Since all $B_{q,m}$ are morphisms of
$U_q\mathfrak{su}_{n,1}$-modules, we need  to verify that $B_{q,m}(f_0)\in
\mathcal{D}(U_n)_q$. This can be easily deduced from the following equivalent
formulation of Proposition \ref{bf0}:
\begin{equation}\label{bfo}
B_{q,t}(f_0)=
(q^2t;q^2)_n\cdot \sum_{k=0}^\infty f_k q^{2k(n+1)}t^k \end{equation}
where $f_k$, $k\in\mathbb{Z}_{\geq0}$, are the finite functions given by
\begin{equation}\label{fk}
f_k=f_k(y_1,y_2,\ldots,y_n)=\left \{\begin{array}{ccl} 1,\,\, &
y_1=q^{2k}
\\ 0,\,\,& {\rm otherwise}\end{array}\right.
\end{equation}
and the series (\ref{bfo}) is convergent in $\mathcal{D}(U_n)'_q$.
Differentiating with respect to $t$ proves
$B_{q,m}(f_0)\in
\mathcal{D}(U_n)_q$.

The equality (\ref{bfo}) and continuity of the
$U_q\mathfrak{su}_{n,1}$-action in $\mathcal{D}(U_n)'_q$ prove also
statement iii) of Proposition \ref{series}.\hfill$\blacksquare$

\medskip

Proposition \ref{series} implies, in particular, that
$B_{q,t}=\mathrm{id}+o(1)$ when $t\rightarrow 0$. By analogy with the
classical case we call the first term of the asymptotic series (\ref{ex})
the $q$-Laplace-Beltrami operator on the quantum ball:

\begin{equation}\label{lb}
\Delta_{n,q}=\frac{q^{-2n}}{1-q^2}\frac{dB_{q,t}}{dt}\mid_{t=0}.
\end{equation}
Clearly, $\Delta_{n,q}:\mathcal{D}(U_n)_q\rightarrow\mathcal{D}(U_n)_q$ is a
morphism of $U_q\mathfrak{su}_{n,1}$-modules.

There are many evidences that this operator should indeed be treated as a
$q$-analogue of the classical Laplace-Beltrami operator in the unit ball.
The results of the present provide an example of such an evidence.

In conclusion, we make also the following remark. The operator
$\Delta_{n,q}$ appeared as the first term of asymptotic of the $q$-Berezin
transform. It turns out that other terms can be explicitly expressed via
$\Delta_{n,q}$. This will be shown in the next section (subsection
\ref{asym}).

\medskip

\subsection{Radial part of the $q$-Laplace-Beltrami operator}

In the classical case the Laplace-Beltrami operator in the unit ball keeps
invariant the space of smooth radial functions, i.~ e. functions depending
on the radius only. The reason is that the radial functions are precisely
the $S(U(n)\times U(1))$-invariant functions with $S(U(n)\times U(1))\subset
SU(n,1)$ being the isotropy group of the centre of the ball. Thus the
'right' $q$-analogue of the radial functions are functions on the quantum
ball which are
$U_q\mathfrak{s}(\mathfrak{u}_n\times\mathfrak{u}_1)$-invariant where
$U_q\mathfrak{s}(\mathfrak{u}_n\times\mathfrak{u}_1)$ is the $*$-Hopf
subalgebra in $U_q\mathfrak{su}_{n,1}$ generated by $E_i, F_i$,
$i=1,2,\ldots,n-1$, and all $K_j^{\pm1}$'s.

It can be proved that any
$U_q\mathfrak{s}(\mathfrak{u}_n\times\mathfrak{u}_1)$-invariant element in
$\mathcal{P}(\mathbb{C}^n)_q$ is a polynomial in
$z_1z_1^*+z_2z_2^*+\ldots+z_nz_n^*$. The idea of the proof is as follows.
$U_q\mathfrak{s}(\mathfrak{u}_n\times\mathfrak{u}_1)$-invariance of a
polynomial, in particular, implies its $U_q\mathfrak{h}$-invariance, i. e.
$K_j$-invariance for any $j$. Obviously, the latter means that the
polynomial depends on $z_1z_1^*, z_2z_2^*,\ldots, z_nz_n^*$ only. One can
write down without difficulties explicit formulas for the action of the
generators $E_i, F_i$, $i=1,2,\ldots,n-1$, on an arbitrary element which
depends on $z_1z_1^*, z_2z_2^*,\ldots, z_nz_n^*$ and find the invariant
elements.

In many computations it is convenient to use the element
$y_1=1-z_1z_1^*-\ldots-z_nz_n^*$ instead of
$z_1z_1^*+z_2z_2^*+\ldots+z_nz_n^*$ since the former quasicommutes with all
the generators $z_i, z_i^*$ (see (\ref{yi})):
$$z_iy_1=q^{-2}y_1z_i, \qquad z^*_iy_1=q^{2}y_1z^*_i.$$
In the sequel we omit the subscript $1$ in the notation for $y_1$.

Using precisely the same arguments, one can show that any
$U_q\mathfrak{s}(\mathfrak{u}_n\times\mathfrak{u}_1)$-invariant finite
function or distribution on the quantum ball depends on $y$ only.

Recall the notation $f_k$, $k\in\mathbb{Z}_{\geq0}$, from the previous
subsection:
$$f_k(y)=\left \{\begin{array}{ccl} 1,\,\,& y=q^{2k}
\\0,\,\,&{\rm otherwise.}\end{array}\right.$$
These functions constitute a basis in the space of radial finite functions
on the quantum ball. They have the following obvious properties: first,
$f_k\cdot f_l=\delta_{kl}f_k$ with $\delta_{kl}$ being the Kronecker symbol,
second, for any distribution $f(y)$
$$f=\sum_{k=0}^\infty f(q^{2k})\cdot f_k$$
where the series converges in $\mathcal{D}(U_n)'_q$.

Our aim now is to compute the action of the $q$-Laplace-Beltrami operator
$\Delta_{n,q}$ on radial finite functions. We would determine the action
completely if we find $\Delta_{n,q}(f_k)$.

\medskip

\begin{proposition}\label{bfk} The operator
$\Delta_{n,q}$ has the following Jacobi matrix form
$$\Delta_{n,q}(f_k)=\frac{q^2}{(1-q^2)^2}\cdot\left((1-q^{2k+2})f_{k+1}-
(1+q^{-2n}-2q^{2k})f_k+ (q^{-2n}-q^{2k-2})f_{k-1}\right)$$ (we assume
$f_{k}\equiv0$ for $k<0$).
\end{proposition}

\medskip

\noindent The proof of the proposition will be given
in the  subsection \ref{proof}.

Proposition \ref{bfk} allows us to apply the $q$-Laplace-Beltrami operator
to any radial function $f=f(y)$. We call restriction of the operator to the
space of radial functions the radial part of the $q$-Laplace-Beltrami
operator and denote it by $\Delta_{n,q}^{(\mathrm{r})}$. Using the explicit
formula from Proposition \ref{bfk}, one can show that the radial part is
given by the following second-order difference operator

\begin{equation}\label{rp}
\Delta_{n,q}^{(\mathrm{r})}=\frac{q^{-n}y^{n+1}}{(yq^2;q^2)_{n-1}}D
y^{-n+1}(yq;q^2)_{n}D
\end{equation}
with $Df(y)=\frac{f(q^{-1}y)-f(qy)}{q^{-1}y-qy}.$

\medskip

\subsection{$q$-Spherical transform}\label{4.3} In this subsection we
describe eigenfunctions of the operator $\Delta_{n,q}^{(\mathrm{r})}$ and
present an explicit formula for expansion in these functions. The associated
$q$-spherical transform should be viewed as a $q$-analogue of the spherical
transform in the unit ball \cite{Helg}. The eigenfunctions of
$\Delta_{n,q}^{(\mathrm{r})}$ appear to be closely related to certain
one-parameter family of the Al-Salam-Chihara polynomials \cite{KS}, and this
observation simplifies proofs of many statements.

Recall the definition of the basic hypergeometric series $_3 \phi_2$
\cite{GR}:

$$_3 \phi_2
\left(\genfrac{}{}{0pt}{}{a_1,\quad a_2,\quad a_3}{b_1,\quad b_2}\quad
q,\quad z\right)=\sum_{n=0}^\infty \frac{(a_1;q)_n \cdot(a_2;q)_n \cdot
(a_3;q)_n}{(b_1;q)_n \cdot(b_2;q)_n\cdot(q;q)_n} z^n.$$ Recall also the
notation $h=\log q^{-2}$ (see Introduction). We define the element
$\phi_\rho(y)\in\mathcal{D}(U_n)'_q$ as follows

\begin{equation}\label{sph}
\phi_\rho(y)=\, _3\phi_2 \left(\genfrac{}{}{0pt}{}{y^{-1},\quad
q^{n+i\rho},\quad q^{n-i\rho}}{q^{2n},\quad 0}\quad q^2,\quad
q^2\right),\qquad \rho\in[0;\frac{2\pi}{h}].
\end{equation}
Obviously, $\phi_\rho(1)=1$. It is a $q$-analogue of the spherical function
in the unit ball as one can see from the following

\medskip

\begin{proposition}\label{eig}

The distribution $\phi_\rho(y)$ is an eigenvector of the operator
$\Delta_{n,q}^{(\mathrm{r})}$:

$$\Delta_{n,q}^{(\mathrm{r})}(\phi_\rho(y))=\lambda(\rho)\cdot\phi_\rho(y)$$
with
$\lambda(\rho)=-q^{2-2n}\frac{(1-q^{n+i\rho})(1-q^{n-i\rho})}{(1-q^{2})^2}$,
$\rho\in[0;\frac{2\pi}{h}]$.
\end{proposition}

\medskip

\noindent{\bf Proof.} Recall the notation $Q_m(x;a,b|q)$ for the
Al-Salam-Chihara polynomials \cite[Section 3.8]{KS}. It is straightforward
that
\begin{equation}\label{al}
\phi_\rho(q^{2m})=\frac{q^{nm}}{(q^{2n};q^2)_m}\cdot
Q_m(\cos\frac{h\rho}{2};q^n,q^n|q^2).
\end{equation}
The statement of the Proposition \ref{eig} is just another formulation of the
recurrence relation for the Al-Salam-Chihara polynomials \cite[(3.8.4)]{KS}.
\hfill$\blacksquare$

\medskip

Let us compute restriction of the invariant integral in the quantum ball
onto the space of radial finite functions.

$$\int_{U_n}f(y)d\nu_q=(q^2;q^2)_n\cdot\sum_{\underline{k}\in P(n)
}f(q^{2k_1})\cdot q^{-2nk_1}\cdot q^{2k_2}\cdot\ldots\cdot
q^{2k_n}=$$$$=(q^2;q^2)_n\cdot\sum_{k=0}^\infty f(q^{2k})\cdot
q^{-2nk}\sum_{0\leq k_n\leq\ldots\leq k_2\leq k}q^{2k_2}\cdot\ldots\cdot
q^{2k_n}.$$ Using the formula (\ref{formula}) and the Jackson integral
(\ref{jackson}), we finally get

\begin{equation}\label{ri}
\int_{U_n}f(y)d\nu_q=
\frac{1-q^{2n}}{1-q^2}\int_0^1f(y)y^{-n-1}(yq^2;q^2)_{n-1}d_{q^2}y.
\end{equation}

Let us denote by $\mathcal{L}_{n,q}$ and $\mathcal{L}'_{n,q}$ the spaces of
radial finite functions and radial distributions on the quantum ball,
respectively. Elements of these spaces can be treated as functions on the
geometric progression $q^{2\mathbb{Z}_{\geq0}}$. We also impose the notation
$\mathcal{L}^2_{n,q}$ for the Hilbert space of 'square integrable' radial
distributions:

$$\mathcal{L}^2_{n,q}=\{f(y)\in\mathcal{L}'_{n,q}\;|
\; \Vert f\Vert_{\mathcal{L}^2}^2=
\frac{1-q^{2n}}{1-q^2}\int_0^1|f(y)|^2y^{-n-1}(yq^2;q^2)_{n-1}
d_{q^2}y<\infty\}.$$

We define the $q$-spherical transform as the map
$\mathscr{F}:\mathcal{L}_{n,q}\rightarrow C^\infty(0;\frac{2\pi}{h})$
given by
\begin{equation}\label{sp}
f(y)\mapsto \mathscr{F}f(\rho)=\frac{1-q^{2n}}{1-q^2}\int_0^1f(y)
\phi_\rho(y)y^{-n-1} (yq^2;q^2)_{n-1}d_{q^2}y.
\end{equation}
It is clear (see the proof of Proposition \ref{eig}) that
$$ \mathscr{F}{f}_k(\rho)=(1-q^{2n})\cdot
q^{-2kn}\cdot(q^{2k+2};q^2)_{n-1}\cdot\phi_\rho(q^{2k})=$$
\begin{equation}\label{Ff}
=(1-q^{2n})\cdot q^{-kn}\cdot
\frac{(q^{2k+2};q^2)_{n-1}}{(q^{2n};q^2)_{k}}\cdot
Q_k(\cos\frac{h\rho}{2};q^n,q^n|q^2).
\end{equation}
 Hence the image of $\mathscr{F}$ is the space of polynomials in
$\cos\frac{h\rho}{2}$.

The following proposition is can be derived from the
spectral decomposition of the operator  $\Delta_{n,q}^{(\mathrm{r})}$,
and a special case of a general result of
\cite[Section 5]{KoeSt}.

\medskip

\begin{proposition}\label{in}

i) The operator $\Delta_{n,q}^{(\mathrm{r})}$ on $\mathcal{L}_{n,q}$ can be
extended to a bounded self-adjoint operator on $\mathcal{L}^2_{n,q}$. It has
simple purely continuous spectrum which coincides with the segment
$[\lambda(\frac{2\pi}{h});\lambda(0)]$ (with $\lambda$ being defined in
Proposition \ref{eig}).

ii) For any finite function $f(y)$

\begin{equation}\label{invers}
f(y)=\frac1{4\pi}\cdot\frac{h}{1-q^{2n}}\cdot\int\limits_0^{\frac{2\pi}{h}}
\mathscr{F}{f}(\rho)\phi_\rho(y)\frac{d\rho}{|c(\rho)|^2}
\end{equation}
where $c(\rho)$ (a $q$-analogue of the Harish-Chandra function) is given by
$$c(\rho)=\frac{\Gamma_{q^2}(n)\Gamma_{q^2}(i\rho)}{\Gamma_{q^2}^2
(\frac{n}{2}+\frac{i\rho}{2})}.$$

iii)  The $q$-spherical transform $\mathscr{F}:\mathcal{L}_{n,q}\rightarrow
C^\infty(0;\frac{2\pi}{h})$ can be extended to a unitary linear operator
$\mathscr{F}:\mathcal{L}^2_{n,q}\rightarrow L^2(\frac{d\rho}{|c(\rho)|^2})$:

\begin{equation}\label{plan}
\frac{1-q^{2n}}
{1-q^2}\int_0^1|f(y)|^2y^{-n-1}(yq^2;q^2)_{n-1}d_{q^2}y=\frac1{4\pi}
\cdot\frac{h}{1-q^{2n}}\cdot \int\limits_0^{\frac{2\pi}{h}}|
\mathscr{F}{f}(\rho)|^2\frac{d\rho}{|c(\rho)|^2}
\end{equation}
(the Plancherel formula).
\end{proposition}

\medskip

Note that, due to (\ref{Ff}), statement iii) can be rewritten as the
orthogonality relations for the Al-Salam-Chihara polynomials \cite[Section
3.8]{KS}:
\begin{equation}\label{ort}
\frac1{4\pi}
\int\limits_0^{\frac{2\pi}{h}}Q_k(\cos\frac{h\rho}{2};q^n,q^n|q^2)
Q_m(\cos\frac{h\rho}{2};q^n,q^n|q^2)\frac{d\rho}{|c(\rho)|^2}=\delta_{km}\cdot
\frac{(q^{2n};q^2)^2_{k}}{h(q^{2k+2};q^2)_{n-1}}.
\end{equation}

\medskip

\subsection{Proof of Proposition \ref{bfk}}\label{proof}
It suffices, due to (\ref{lb}), to verify that
\begin{equation}\label{modulo-t2}
B_{q,t}(f_k)=f_k+t\frac{q^{2n+2}}{(1-q^2)}\cdot\left((1-q^{2k+2})f_{k+1}-
(1+q^{-2n}-2q^{2k})f_k+ (q^{-2n}-q^{2k-2})f_{k-1}\right)
\end{equation}
modulo $t^2$. For that purpose we compute first $T_{f_{k}}$.

Recall the notation $\mathbb{C}[\mathbb{C}^n]_{q,\alpha}$
(subsection \ref{2.5}).
Define a $\mathbb{Z}_{\geq0}$-grading in the space
$\mathbb{C}[\mathbb{C}^n]_{q,\alpha}$ as follows

$$\mathbb{C}[\mathbb{C}^n]^{(m)}_{q,\alpha}=\mathrm{linear}\;\mathrm{span}\;
\mathrm{of}\;\mathbf{z}^{\underline{m}},\;|\underline{m}|=m.$$ We need first

\medskip

\begin{lemma}\label{tfk}
$$T_{f_k}\mid_{\mathbb{C}[\mathbb{C}^n]^{(m)}_{q,\alpha}}=q^{2(k-m)(\alpha+1)}
\cdot
\frac{(q^{2\alpha+2};q^2)_{n+m}\cdot(q^{2k-2m+2};q^2)_{n+m-1}}
{(q^2;q^2)_{n+m-1}}.$$
\end{lemma}

\medskip

\noindent{\bf Proof of the lemma.} Using the same arguments as in the
classical case, it can be proved that
$\mathbb{C}[\mathbb{C}^n]^{(m)}_{q,\alpha}$, $m\in\mathbb{Z}_{\geq0}$, are
pairwise non-isomorphic irreducible
$U_q\mathfrak{s}(\mathfrak{u}_n\times\mathfrak{u}_1)$-modules (we mean
restriction of the representation $\pi_\alpha$ onto
$U_q\mathfrak{s}(\mathfrak{u}_n\times\mathfrak{u}_1)$). Due to the
$U_q\mathfrak{s}(\mathfrak{u}_n\times\mathfrak{u}_1)$-invariance of $f_k$
(with respect to the untwisted
$U_q\mathfrak{s}(\mathfrak{u}_n\times\mathfrak{u}_1)$-action), Proposition
\ref{toepinv}, and the Schur lemma\footnote{Since
$\mathbb{C}[\mathbb{C}^n]^{(m)}_{q,\alpha}$ is an irreducible
$U_q\mathfrak{s}(\mathfrak{u}_n\times\mathfrak{u}_1)$-module, the algebra
homomorphism
$U_q\mathfrak{s}(\mathfrak{u}_n\times\mathfrak{u}_1)\to\mathrm{End}
(\mathbb{C}[\mathbb{C}^n]^{(m)}_{q,\alpha})$ is surjective, however the only
operators commuting with the full matrix algebra are constants.}

$$ T_{f_k}\mid_{\mathbb{C}[\mathbb{C}^n]^{(m)}_{q,\alpha}}=c_k^m.$$ To compute
the constant $c_k^m$ we apply $T_{f_k}$ to a distinguished vector in
$\mathbb{C}[\mathbb{C}^n]^{(m)}_{q,\alpha}$, for example, $z_n^m$:

$$c_k^m=\frac{(T_{f_k}z_n^m,z_n^m)_{\alpha,q}}{(z_n^m,z_n^m)_{\alpha,q}}=
\frac{(f_kz_n^m,z_n^m)_{\alpha,q}}{(z_n^m,z_n^m)_{\alpha,q}}.$$ By Proposition
\ref{norms} the denominator is equal to
$\frac{(q^2;q^2)_m}{(q^{2n+2\alpha+2};q^2)_m}$. What remains is to compute the
numerator:

$$(T_{f_k}z_n^m,z_n^m)_{\alpha,q}=\int_{U_n}z_n^{*m}f_kz_n^mdm_{\alpha,q}=
\int_{U_n}f_{k-m}(q^2y_n;q^2)_mdm_{\alpha,q}=$$
$$=(q^{2\alpha+2};q^2)_n\sum_{\underline{k}\in
P(n)}f_{k-m}(q^{2k_1})(q^{2+2k_n};q^2)_mq^{2k_1(\alpha+1)}q^{2k_2}\ldots
q^{2k_n}=$$$$=(q^{2\alpha+2};q^2)_nq^{2(k-m)(\alpha+1)} \sum_{0\leq
k_n\ldots\leq k_2 \leq k-m}(q^{2+2k_n};q^2)_mq^{2k_2}\ldots q^{2k_n}.$$ To
continue computation we need the following simple formula which can be
proved by induction:
\begin{equation}\label{formula}
\sum_{a\leq l_n\ldots\leq l_1 \leq b}q^{2l_1}\ldots
q^{2l_n}=q^{2an}\frac{(q^{2b-2a+2};q^2)_n}{(q^2;q^2)_n}.
\end{equation}
By this formula $(q^{2+2k_n};q^2)_m=(q^2;q^2)_m
\sum_{0\leq l_m\ldots\leq l_1 \leq
k_n}q^{2l_1}\ldots q^{2l_m}$. Applying (\ref{formula}) one more time,
we finally get

\begin{equation*}
(T_{f_k}z_n^m,z_n^m)_{\alpha,q}=(q^{2\alpha+2};q^2)_nq^{2(k-m)(\alpha+1)}
(q^2;q^2)_m
\frac{(q^{2k-2m+2};q^2)_{n+m-1}}{(q^2;q^2)_{n+m-1}}.
\end{equation*}
\hfill$\blacksquare$

\medskip

Let $P_m$ denotes the orthogonal projection in
$\mathbb{C}[\mathbb{C}^n]_{q,\alpha}$
onto $\mathbb{C}[\mathbb{C}^n]^{(m)}_{q,\alpha}$. The above lemma says that
\begin{equation}\label{formu}
T_{f_{k}}=\sum_{m=0}^k q^{2(k-m)(\alpha+1)}\cdot
\frac{(q^{2\alpha+2};q^2)_{n+m}
\cdot
(q^{2k-2m+2};q^2)_{n+m-1}}{(q^2;q^2)_{n+m-1}} \cdot P_m.
\end{equation}

\medskip

\begin{lemma}\label{csp}
The covariant symbol $\sigma(P_m)$ is given by

$$\sigma(P_m)=q^{-2m(\alpha+n+1)}\cdot
\frac{(q^{2\alpha+2n+2};q^2)_{m}}{(q^2;q^2)_{m}} \cdot
y^{\alpha+n+1}\cdot(yq^{-2m+2};q^2)_m.$$
\end{lemma}

\medskip

\noindent{\bf Proof of the lemma.} Since
$\mathbb{C}[\mathbb{C}^n]^{(m)}_{q,\alpha}$, $m\in\mathbb{Z}_{\geq0}$, are
pairwise non-isomorphic irreducible
$U_q\mathfrak{s}(\mathfrak{u}_n\times\mathfrak{u}_1)$-modules, the
projections $P_m$ are $U_q\mathfrak{s}(\mathfrak{u}_n\times
\mathfrak{u}_1)$-module morphisms. Thus by Proposition \ref{invcov}
$\sigma(P_m)$ should be a function of $y$. Denote it by $p_m(y)$. Recall
that $p_m(y)=\sum_{k=0}^\infty p_m(q^{2k})f_k(y)$. The coefficients
$p_m(q^{2k})$ can be derived from the equalities

$$\int_{U_n}p_m(y)\cdot f_ld\nu_q=\mathrm{Tr}_q(P_m\cdot T_{f_l}),\quad
l\in\mathbb{Z}_{\geq0}.$$ Using the property $f_k\cdot f_l=\delta_{kl}f_k$
and lemma \ref{tfk}, we can rewrite the latter equality as follows:

$$p_m(q^{2l})\cdot (q^2;q^2)_n\cdot q^{-2nl}\cdot
\sum_{0\leq k_n\ldots\leq k_2
\leq
l}q^{2k_2}\ldots q^{2k_n}=$$$$=q^{2(l-m)(\alpha+1)}\cdot
\frac{(q^{2\alpha+2};q^2)_{n+m}\cdot(q^{2l-2m+2};q^2)_{n+m-1}}
{(q^2;q^2)_{n+m-1}}\cdot\mathrm{Tr}_q(P_m).$$ By (\ref{qtra}) $
\mathrm{Tr}_q(P_m)=\frac{(q^2;q^2)_{n}}{(q^{2\alpha+2};q^2)_{n}}\cdot
q^{-2nm}\cdot\sum_{0\leq k_n\ldots\leq k_2 \leq m}q^{2k_2}\ldots q^{2k_n}$.
Using (\ref{formula}) we obtain $
p_m(q^{2l})=
\frac{(q^{2\alpha+2n+2};q^2)_{m}\cdot(q^{2l-2m+2};q^2)_m}{(q^2;q^2)_{m}}
\cdot q^{2(l-m)(\alpha+n+1)}$. \hfill$\blacksquare$

\medskip

We continue the proof of our proposition. Lemmas \ref{tfk} and \ref{csp} now
imply
\begin{equation}\label{for}
\begin{split}
B_{q,\alpha}(f_k)&=q^{2k(\alpha+1)}\cdot
y^{\alpha+n+1}\cdot\sum_{m=0}^k q^{-4m(\alpha+1)-2mn}\cdot
\\
&\quad \cdot\frac{(q^{2\alpha+2};q^2)_{n+m}(q^{2k-2m+2};q^2)_{n+m-1}
(q^{2\alpha+2n+2};q^2)_{m}}{(q^2;q^2)_{n+m-1}(q^2;q^2)_{m}}
\cdot(yq^{-2m+2};q^2)_m.
\end{split}\end{equation}
Note that $(yq^{-2m+2};q^2)_m=\sum_{l=0}^\infty (q^{2l-2m+2};q^2)_m\cdot
f_l=\sum_{l=m}^\infty (q^{2l-2m+2};q^2)_m\cdot f_l$. Consequently
\begin{equation*}
\begin{split}
&\qquad y^{\alpha+n+1}\cdot(yq^{-2m+2};q^2)_m=\sum_{l=m}^\infty
q^{2l(\alpha+n+1)}\cdot(q^{2l-2m+2};q^2)_m\cdot f_l=\\ &=\sum_{r=0}^\infty
q^{2(m+r)(\alpha+n+1)}\cdot(q^{2r+2};q^2)_m\cdot f_{m+r},\end{split}
\end{equation*} and
\begin{equation*}
\begin{split}
B_{q,\alpha}(f_k)&=\sum_{m=0}^k\sum_{r=0}^\infty q^{2\alpha(k-m+r)}\cdot
q^{2k-2m+2rn+2r}\cdot\\ &\qquad
\cdot\frac{(q^{2\alpha+2};q^2)_{n+m}(q^{2k-2m+2};q^2)_{n+m-1}
(q^{2\alpha+2n+2};q^2)_{m}(q^{2r+2};q^2)_m}{(q^2;q^2)_{n+m-1}(q^2;q^2)_{m}}
\cdot
f_{m+r}.
\end{split}
\end{equation*}
 The substitution $t=q^{2\alpha}$ gives
\begin{equation*}
\begin{split}
B_{q,t}(f_k)&=
\sum_{m=0}^k\sum_{r=0}^\infty t^{k-m+r}\cdot q^{2k-2m+2rn+2r}\cdot\\
&\cdot\frac{(tq^{2};q^2)_{n+m}(q^{2k-2m+2};q^2)_{n+m-1}
(tq^{2n+2};q^2)_{m}(q^{2r+2};q^2)_m}
{(q^2;q^2)_{n+m-1}(q^2;q^2)_{m}}\cdot f_{m+r}.
\end{split}
\end{equation*}
To finish the proof of (\ref{modulo-t2}) we have to compute the coefficient
at $t$ in the last series. Due to presence of the multiply $t^{k-m+r}$ in
the terms, the only terms which might give a non-zero contribution to this
coefficient correspond to the values $(m,r)=(k,0)$, $(m,r)=(k,1)$, or
$(m,r)=(k-1,0)$ which we treat separately

1)$(m,r)=(k,0)$: the corresponding term is equal to

$$(tq^{2};q^2)_{n+k}(tq^{2n+2};q^2)_{k}f_k=\left(1-t\frac{q^{2n+2}}{(1-q^2)}
(1+q^{-2n}-2q^{2k})\right)f_k+o(t).$$

2)$(m,r)=(k,1)$: the corresponding term is equal to

$$t\cdot q^{2n+2}\cdot
(tq^{2};q^2)_{n+k}(tq^{2n+2};q^2)_{k}\frac{1-q^{2k+2}}{1-q^2}f_{k+1}=
t\cdot\frac{q^{2n+2}}{(1-q^2)}\cdot(1-q^{2k+2})f_{k+1}+o(t).$$

3)$(m,r)=(k-1,0)$: the corresponding term is equal to

$$t\cdot q^{2}\cdot
(tq^{2};q^2)_{n+k-1}(tq^{2n+2};q^2)_{k-1}\frac{1-q^{2n+2k-2}}{1-q^2}f_{k-1}=
t\cdot\frac{q^{2n+2}}{(1-q^2)}\cdot(q^{-2n}-q^{2k-2})f_{k-1}+o(t).$$ This
finishes the proof.

\bigskip

\section{Further properties of the $q$-Berezin transform}\label{5}

In this section we study further the $q$-Berezin transform. Namely, we
consider its restriction $B^{(\mathrm{r})}_{q,\alpha}$ onto the space of
radial functions. We prove that $B^{(\mathrm{r})}_{q,\alpha}$ is extended to
a bounded selfadjoint operator on $\mathcal{L}^2_{n,q}$ which commutes with
the radial part $\Delta^{(\mathrm{r})}_{n,q}$ of the $q$-Laplace-Beltrami
operator. Since the latter has a simple spectrum,
$B^{(\mathrm{r})}_{q,\alpha}$ is a function of
$\Delta^{(\mathrm{r})}_{n,q}$. We find the function explicitly. We also
present an asymptotic expansion for the $q$-Berezin transform at the limit
$t=q^{2\alpha}\to0$ mentioned at the end of subsection \ref{4.1}.

\medskip

\subsection{Boundedness of the $q$-Berezin transform}

Let $B^{(\mathrm{r})}_{q,\alpha}:\mathcal{L}_{n,q}\rightarrow
\mathcal{L}'_{n,q}$ be restriction of the $q$-Berezin transform onto the
space $\mathcal{L}_{n,q}$ of finite radial functions on the quantum ball.

\medskip

\begin{proposition}\label{sy}
$B^{(\mathrm{r})}_{q,\alpha}$ can be extended to a bounded self-adjoint
operator on $\mathcal{L}^2_{n,q}$. It is a function of
$\Delta^{(\mathrm{r})}_{n,q}$. The operator $\mathscr{F}\cdot
B^{(\mathrm{r})}_{q,\alpha}\cdot\mathscr{F}^{-1}$ on
$L^2(\frac{d\rho}{|c(\rho)|^2})$ is the multiplication by the (bounded)
function
\begin{equation}\label{bqa}
b_{q,\alpha}(\rho)=
\frac{(q^{2+2\alpha};q^2)_\infty\cdot(q^{2n+2+2\alpha};q^2)_\infty}
{(q^{n+2+2\alpha+i\rho};q^2)_\infty\cdot(q^{n+2+2\alpha-i\rho};q^2)_\infty}.
\end{equation}
\end{proposition}

\medskip

\noindent{\bf Proof.} We divide the proof into three lemmas.

\medskip

\begin{lemma}
$B^{(\mathrm{r})}_{q,\alpha}(\mathcal{L}_{n,q})\subset \mathcal{L}^2_{n,q}$
and one has the equality
\begin{equation}\label{com}
\Delta^{(\mathrm{r})}_{n,q}\cdot B^{(\mathrm{r})}_{q,\alpha}
=B^{(\mathrm{r})}_{q,\alpha}\cdot\Delta^{(\mathrm{r})}_{n,q}
\end{equation}
of linear maps from $\mathcal{L}_{n,q}$ to $\mathcal{L}^2_{n,q}$.
\end{lemma}

\medskip

\noindent{\bf Sketch of a proof.} The inclusion is due to the formula
(\ref{for}). Indeed, one has to show that
$y^{\alpha+n+1}\cdot(yq^{-2m+2};q^2)_m\in\mathcal{L}^2_{n,q}$ for any $m$,
and this is clear be the definition.

To prove the equality (\ref{com}), we have to show that
$\Delta^{(\mathrm{r})}_{n,q}\cdot B^{(\mathrm{r})}_{q,\alpha}(f_k)
=B^{(\mathrm{r})}_{q,\alpha}\cdot\Delta^{(\mathrm{r})}_{n,q}(f_k)$ for any
$k$. The left-hand side may computed by successive application of
(\ref{for}) and (\ref{rp}) while the right-hand one may be computed via
Proposition \ref{bfk} and (\ref{for}). \hfill$\blacksquare$

\medskip

\begin{lemma} Suppose $B$ is a linear operator from
$\mathcal{L}_{n,q}$ to $\mathcal{L}^2_{n,q}$ which satisfies the properties
$$\Delta^{(\mathrm{r})}_{n,q}\cdot B =B\cdot\Delta^{(\mathrm{r})}_{n,q};$$
$$Bf_0=
(q^{2\alpha+2};q^2)_n\cdot y^{\alpha+n+1}.$$ Then
$B=B^{(\mathrm{r})}_{q,\alpha}.$
\end{lemma}

\medskip

\noindent{\bf Proof.} The statement is a simple consequence of the equality

$$\mathcal{L}_{n,q}=\text{linear span of}\,\,
\{(\Delta^{(\mathrm{r})}_{n,q})^mf_0\}_{m\in\mathbb{Z}_{\geq0}}$$ which in
turn may be deduced easily from Proposition \ref{bfk}. \hfill$\blacksquare$

\medskip

\begin{lemma}
The bounded selfadjoint operator $\mathscr{F}^{-1}\cdot
b_{q,\alpha}(\rho)\cdot\mathscr{F}$ ($b_{q,\alpha}$ is given by (\ref{bqa}))
on $\mathcal{L}^2_{n,q}$ possesses the properties from the previous lemma.
\end{lemma}

\medskip

\noindent{\bf Proof.}The first property holds trivially. Let us prove the
second one.

Let $(\cdot \;,\;\cdot)_{L^2}$, $(\cdot \;,\;\cdot)_{\mathcal{L}^2}$ be the
inner products in $L^2(\frac{d\rho}{|c(\rho)|^2})$ and
$\mathcal{L}^2_{n,q}$, respectively. It is sufficient to show that

$$(b_{q,\alpha}(\rho)\cdot\mathscr{F}f_0
\;,\;\mathscr{F}f_k)_{L^2}=(q^{2\alpha+2};q^2)_n\cdot (y^{\alpha+n+1}
\;,\;f_k)_{\mathcal{L}^2}$$ for any $k\in\mathbb{Z}_{\geq0}$. Recall
(\ref{Ff}) that $\mathscr{F}{f}_k(\rho)=(1-q^{2n})\cdot q^{-kn}\cdot
\frac{(q^{2k+2};q^2)_{n-1}}{(q^{2n};q^2)_{k}}\cdot Q_k(\rho)$ with
$Q_k(\rho):=Q_k(\cos\frac{h\rho}{2};q^n,q^n|q^2)$. We rewrite the condition
as follows
\begin{equation}\label{bqal}
(b_{q,\alpha}(\rho) \;,\;Q_k(\rho))_{L^2}=
q^{k(n+2+2\alpha)}\cdot\frac{(q^{2\alpha+2};q^2)_{n}}{(q^{2};q^2)_{n}}.
\end{equation}
Remind the following formula for the generating function of the
Al-Salam-Chihara polynomials \cite[formula (3.8.13)]{KS}:

$$\sum_{j=0}^\infty \frac{z^j}{(q^2;q^2)_{j}}\cdot
Q_j(\rho)=\frac{(q^nz;q^2)_\infty\cdot(q^{n}z;q^2)_\infty}
{(q^{i\rho}z;q^2)_\infty\cdot(q^{-i\rho}z;q^2)_\infty}.$$ Substitution
$z:=q^{2\alpha+2+n}$ gives $(q^{2\alpha+2};q^2)_n\sum_{j=0}^\infty
\frac{q^{j(2\alpha+n+2)}}{(q^2;q^2)_{j}}\cdot Q_j(\rho)=b_{q,\alpha}(\rho)$.
Now validity of (\ref{bqal}) follows from the latter equality and the
orthogonality relations (\ref{ort}) for the Al-Salam-Chihara
polynomials.\hfill$\blacksquare$

\medskip

Proposition \ref{sy} follows directly from the last two lemmas.

\medskip

\subsection{Asymptotic expansion of the $q$-Berezin
transform}\label{asym}

\begin{proposition}\label{as}
  For any $f\in\mathcal{D}(U_n)_q$
$$B_{q,t}(f)=(q^2t;q^2)_n\cdot\sum_{j=0}^\infty t^j\cdot q^{2j}
\cdot\frac{(q^{2j+2};q^2)_{n-1}}{(q^2;q^2)_{n-1}}\cdot p_j(\Delta_{n,q})f$$
with
$$p_j(\Delta_{n,q})= \sum_{l=0}^j \frac{(q^{-2j};q^2)_l\cdot
q^{2l}}{(q^{2n};q^2)_l\cdot (q^2;q^2)_l} \cdot
\prod_{m=0}^{l-1}\left((1-q^{2m})(1-q^{2m+2n})-
q^{2m+2n-2}(1-q^2)^2\Delta_{n,q}\right).$$
\end{proposition}

\medskip

\noindent{\bf Proof.} Recall (see (\ref{bf0})) that

$$B_{q,\alpha}(f_0)=(q^{2\alpha+2};q^2)_n\cdot
y^{\alpha+n+1}=(q^{2\alpha+2};q^2)_n\cdot\sum_{j=0}^\infty
q^{2j(\alpha+n+1)}\cdot f_j$$ or, in terms of $t=q^{2\alpha}$,
$$B_{q,t}(f_0)=(tq^{2};q^2)_n\cdot \sum_{j=0}^\infty t^{j}q^{2j(n+1)}\cdot
f_j.$$ Since any morphism of $U_q\mathfrak{su}_{n,1}$-modules is determined
uniquely by its value on the vector $f_0$, it suffices to prove that $
p_j(\Delta_{n,q})f_0=q^{2jn}\cdot
\frac{(q^2;q^2)_{n-1}}{(q^{2j+2};q^2)_{n-1}}\cdot f_j$ or, equivalently, $$
p_j(\Delta^{(\mathrm{r})}_{n,q})f_0=q^{2jn}\cdot
\frac{(q^2;q^2)_{n-1}}{(q^{2j+2};q^2)_{n-1}}\cdot f_j.$$ Let us apply the
$q$-spherical transform $\mathscr{F}$ to the both hand-sides of the latter
equality. We get
\begin{equation}\label{pj}
 p_j(\lambda(\rho))\mathscr{F}{f}_0(\rho)=q^{2jn}\cdot
\frac{(q^2;q^2)_{n-1}}{(q^{2j+2};q^2)_{n-1}}\cdot \mathscr{F}{f}_j(\rho).
\end{equation}
Recall (see (\ref{Ff})) that $ \mathscr{F}{f}_k(\rho)=(1-q^{2n})\cdot
q^{-2kn}\cdot(q^{2k+2};q^2)_{n-1}\cdot\phi_\rho(q^{2k})$. This formula
reduces proving (\ref{pj}) to proving the equality
\begin{equation}\label{pphi}
p_j(\lambda(\rho))=\phi_\rho(q^{2j}), \end{equation} and the latter is just
a straightforward computation.\hfill$\blacksquare$

\bigskip

\section{An application: orthogonality relations for continuous dual $q$-Hahn
polynomials}\label{6}

The aim of this section is to describe one application of our results to the
theory of basic orthogonal polynomials. Namely, we use the $q$-Berezin and
the $q$-spherical transform to obtain orthogonality relations for certain
two-parameter family of the so-called continuous dual $q$-Hahn polynomials
(see \cite{KS}). Of course, this result is not new. However we believe that
our approach might be interesting.

Throughout this section $\alpha$ is a fixed number.

\medskip

\subsection{Auxiliary results} The aim of this subsection is to derive some
useful consequences of Proposition \ref{sy}.

\medskip

\begin{proposition}\label{invert}
$B^{(\mathrm{r})}_{q,\alpha}$ is extended to an invertible operator on
$\mathcal{L}^2_{n,q}$.
\end{proposition}

\medskip

\noindent{\bf Proof.} It suffices to observe that the 'symbol'
$b_{q,\alpha}(\rho)$ is invertible:

$$b_{q,\alpha}(\rho)\geq
\frac{(q^{2+2\alpha};q^2)_\infty\cdot(q^{2n+2+2\alpha};q^2)_\infty}
{(-q^{n+2+2\alpha};q^2)_\infty\cdot(-q^{n+2+2\alpha};q^2)_\infty}.$$
\hfill$\blacksquare$

\medskip

To go further, we introduce some auxiliary notations. Denote by
$\mathcal{L}_{\mathrm{Op}}$ and $\mathcal{L}'_{\mathrm{Op}}$ the subspaces
of $U_q\mathfrak{s}(\mathfrak{u}_n\times\mathfrak{u}_1)$-module morphisms in
$\mathrm{End}_0(\mathbb{C}[\mathbb{C}^n]_{q,\alpha})$ and
$\mathrm{End}(\mathbb{C}[\mathbb{C}^n]_{q,\alpha})$, respectively. Since the
subspaces $\mathbb{C}[\mathbb{C}^n]^{(m)}_{q,\alpha}$,
$m\in\mathbb{Z}_{\geq0}$, in $\mathbb{C}[\mathbb{C}^n]_{q,\alpha}$ are
pairwise non-isomorphic irreducible $U_q\mathfrak{s}(\mathfrak{u}_n\times
\mathfrak{u}_1)$-modules (see the proof of lemma \ref{tfk}), the vector
space $\mathcal{L}_{\mathrm{Op}}$ is generated by the orthogonal projections
$P_m$ onto the subspaces $\mathbb{C}[\mathbb{C}^n]^{(m)}_{q,\alpha}$.
$\mathcal{L}'_{\mathrm{Op}}$ is the space of infinite series of the form

$$T=\sum_ma_mP_m, \quad a_m\in \mathbb{C}.$$ Elements of
$\mathcal{L}_{\mathrm{Op}}$
and $\mathcal{L}'_{\mathrm{Op}}$ play the role of 'radial' elements in
$\mathrm{End}_0(\mathbb{C}[\mathbb{C}^n]_{q,\alpha})$ and
$\mathrm{End}(\mathbb{C}[\mathbb{C}^n]_{q,\alpha})$.

Let $\mathcal{L}^2_{\mathrm{Op}}$ be the subspace in
$\mathcal{L}'_{\mathrm{Op}}$ of '$q$-Hilbert-Schmidt' operators:

$$
\mathcal{L}^2_{\mathrm{Op}}=
\{T\in\mathcal{L}'_{\mathrm{Op}}\;|\;\mathrm{Tr}_q(T^*\cdot
T)<\infty\}.$$ Note that the triple ($\mathcal{L}_{\mathrm{Op}}$,
$\mathcal{L}'_{\mathrm{Op}}$, $\mathcal{L}^2_{\mathrm{Op}})$ is very similar
to the triple $(\mathcal{L}_{n,q}$, $\mathcal{L}'_{n,q}$,
$\mathcal{L}^2_{n,q})$ introduced in subsection \ref{4.3}.

Recall (subsection \ref{3.3}) the notation $\sigma$ for the linear map from
$\mathrm{End}_0(\mathbb{C}[\mathbb{C}^n]_{q,\alpha})$ to
$\mathcal{D}(U_n)'_q$ which sends endomorphisms to their covariant symbols.
Due to Proposition \ref{invcov} $\sigma(\mathcal{L}_{\mathrm{Op}})\subset
\mathcal{L}'_{n,q}$. Let us denote restriction of $\sigma$ onto
$\mathcal{L}_{\mathrm{Op}}$ by $\sigma^{(\mathrm{r})}$.

\medskip

\begin{proposition}\label{boundsigma}
The image of operator $\sigma^{(\mathrm{r})}$ lies in $\mathcal{L}^2_{n,q}$.
Moreover, $\sigma^{(\mathrm{r})}$ can be extended to a bounded invertible
operator from $\mathcal{L}^2_{\mathrm{Op}}$ to $\mathcal{L}^2_{n,q}$.
\end{proposition}

\medskip

\noindent{\bf Proof.} Recall (subsection \ref{3.3}) that the map $\sigma$ is
defined via the equality

$$\int_{U_n}\sigma(T)\cdot fd\nu_q=\mathrm{Tr}_q(T\cdot T_f)$$ which should be
fulfilled for any $f\in\mathcal{D}(U_n)_q$. Let $T^{(\mathrm{r})}$ be
restriction onto $\mathcal{L}_{n,q}$ of the operator from
$\mathcal{D}(U_n)_q$ to
$\mathrm{End}_0(\mathbb{C}[\mathbb{C}^n]_{q,\alpha})$ which sends finite
functions to the corresponding Toeplitz operators. By Proposition
\ref{toepinv}
$T^{(\mathrm{r})}(\mathcal{L}_{n,q})\subset\mathcal{L}'_{\mathrm{Op}}$. More
precisely, by lemma \ref{tfk}
\begin{equation}\label{inclu}
  T^{(\mathrm{r})}(\mathcal{L}_{n,q})\subset\mathcal{L}_{\mathrm{Op}}.
\end{equation}
In particular, $T^{(\mathrm{r})}$ can be considered as a densely defined
operator from $\mathcal{L}^2_{n,q}$ into $\mathcal{L}^2_{\mathrm{Op}}.$ The
above definition of $\sigma$ implies the following equivalent definition of
$\sigma^{(\mathrm{r})}$: for an element $T\in\mathcal{L}_{\mathrm{Op}}$ one
has $\sigma^{(\mathrm{r})}(T)=t(y)$ iff
\begin{equation}\label{eqdef}
\frac{1-q^{2n}}{1-q^2}\int_0^1t(y)f(y)y^{-n-1}(yq^2;q^2)_{n-1}d_{q^2}y=
\mathrm{Tr}_q(T\cdot T^{(\mathrm{r})}(f(y)))
\end{equation}
for any $f(y)\in\mathcal{L}_{n,q}.$ Let us prove that $T^{(\mathrm{r})}$ can
be extended up to a bounded operator from $\mathcal{L}^2_{n,q}$ to
$\mathcal{L}^2_{\mathrm{Op}}$. Suppose $f(y)\in\mathcal{L}_{n,q}$. Then, by
(\ref{inclu}) and (\ref{eqdef}) $$
\mathrm{Tr}_q((T^{(\mathrm{r})}(f(y)))^*\cdot T^{(\mathrm{r})}(f(y)))=
\frac{1-q^{2n}}{1-q^2}\int_0^1\sigma^{(\mathrm{r})}\cdot
T^{(\mathrm{r})}(f(y))\overline{f(y)}y^{-n-1}(yq^2;q^2)_{n-1}d_{q^2}y=$$
$$=\frac{1-q^{2n}}{1-q^2}\int_0^1 B^{(\mathrm{r})}_{q,\alpha}(f(y))
\overline{f(y)}y^{-n-1}(yq^2;q^2)_{n-1}d_{q^2}y=
(B^{(\mathrm{r})}_{q,\alpha}f, f)_{\mathcal{L}^2 }
\leq\|B^{(\mathrm{r})}_{q,\alpha}\| (f, f)_{\mathcal{L}^2 }.$$ To prove
boundedness of $\sigma^{(\mathrm{r})}$ it remains to observe that, by
(\ref{eqdef}), $\sigma^{(\mathrm{r})}$ coincides on
$\mathcal{L}_{\mathrm{Op}}$ with $(T^{(\mathrm{r})})^*$. Invertibility of
$\sigma^{(\mathrm{r})}$ may be deduced using similar arguments and
Proposition \ref{invert}. \hfill$\blacksquare$

\medskip

\subsection{Orthogonality relations}

First we describe an idea of producing the orthogonality relations. It is
based on the following statement.

\medskip

\begin{proposition}\label{main}
The operator $$U_{\alpha}:\mathcal{L}^2_{\mathrm{Op}}\rightarrow
L^2(\frac{d\rho}{|c(\rho)|^2}),\qquad T\mapsto
\frac1{\sqrt{b_{q,\alpha}(\rho)}}\cdot \mathscr{F}\sigma^{(\mathrm{r})}(T)$$
is unitary.
\end{proposition}

\medskip

\noindent{\bf Remark.} In the classical case the operator $U_{\alpha}$ (the
product of the spherical transform with the unitary part of the covariant
symbol map) was studied in \cite{jp+gkz-pl3}.

\medskip

\noindent{\bf Proof.} By Proposition \ref{sy}

$$\frac1{b_{q,\alpha}(\rho)}\cdot\mathscr{F}=\mathscr{F}
(B^{(\mathrm{r})}_{q,\alpha})^{-1}.$$ Remind the notation $(\cdot
\;,\;\cdot)_{L^2}$, $(\cdot \;,\;\cdot)_{\mathcal{L}^2}$ for the inner
products in $L^2(\frac{d\rho}{|c(\rho)|^2})$ and $\mathcal{L}_{n,q}$,
respectively. Then
$$
\left(
\frac1{\sqrt{b_{q,\alpha}(\rho)}}\cdot\mathscr{F}\sigma^{(\mathrm{r})}(T)
\;,\;\frac1{\sqrt{b_{q,\alpha}(\rho)}}\cdot
\mathscr{F}\sigma^{(\mathrm{r})}(T)\right)_{L^2}=$$
$$=
\left(
\frac1{b_{q,\alpha}(\rho)}\cdot\mathscr{F}\sigma^{(\mathrm{r})}(T) \;,\;
\mathscr{F}\sigma^{(\mathrm{r})}(T)\right)_{L^2}=
\left(\mathscr{F}(B^{(\mathrm{r})}_{q,\alpha})^{-1}\sigma^{(\mathrm{r})}(T)
\;,\; \mathscr{F}\sigma^{(\mathrm{r})}(T)\right)_{L^2}=$$$$=
\left((B^{(\mathrm{r})}_{q,\alpha})^{-1} \sigma^{(\mathrm{r})}(T)
\;,\;\sigma^{(\mathrm{r})}(T)\right)_{\mathcal{L}^2}.$$ The result then
follows since
$(B^{(\mathrm{r})}_{q,\alpha})^{-1}=(\sigma^{(\mathrm{r})*})^{-1}\cdot
(\sigma^{(\mathrm{r})})^{-1}.$ \hfill$\blacksquare$

\medskip

Recall that the projections $P_m$, $m\in\mathbb{Z}_{\geq0}$ constitute an
orthogonal basis in the Hilbert space $\mathcal{L}^2_{Op}$. Proposition
\ref{main} implies that $\mathscr{F}\sigma^{(\mathrm{r})}(P_m)$,
$m\in\mathbb{Z}_{\geq0}$ constitute an orthogonal basis in
$L^2(\frac{d\rho}{b_{q,\alpha}(\rho)|c(\rho)|^2})$. Our intention now is to
show that $\mathscr{F}\sigma^{(\mathrm{r})}(P_m)$ are very close to certain
continuous dual $q$-Hahn polynomials, and the above observations will give
us the orthogonality relations for these polynomials.

\medskip

\begin{proposition}\label{fspm}Let
$P_m (\rho)=\mathscr{F}\sigma^{(\mathrm{r})}(P_m)(\rho)$. Then
$$P_m(\rho)=\frac{(q^{2};q^2)_n\cdot(q^{2n+2+2\alpha};q^2)^2_\infty}
{(q^{n+2+2\alpha+i\rho};q^2)_\infty
\cdot(q^{n+2+2\alpha-i\rho};q^2)_\infty}\cdot
q^{-2m(n+\alpha)}\cdot \frac{(q^{2n+2+2\alpha};q^2)_m}{(q^{2};q^2)_m}\cdot$$
$$\cdot _3{\phi} _2 \left(\genfrac{}{}{0pt}{}{q^{-2m},\quad
q^{n+2+2\alpha+i\rho},\quad q^{n+2+2\alpha-i\rho}}{q^{2n+2+2\alpha},\quad
q^{2n+2+2\alpha}}\quad q^2,\quad q^2\right).$$
\end{proposition}

\medskip

\noindent{\bf Proof.} Remind (lemma \ref{csp}) that we have already computed
$\sigma^{(\mathrm{r})}(P_m)$

$$\sigma^{(\mathrm{r})}(P_m)=q^{-2m(\alpha+n+1)}\cdot
\frac{(q^{2\alpha+2n+2};q^2)_{m}}{(q^2;q^2)_{m}} \cdot
y^{\alpha+n+1}\cdot(yq^{-2m+2};q^2)_m.$$ By the $q$-binomial formula
\cite{GR} the equality can be rewritten as follows

$$\sigma^{(\mathrm{r})}(P_m)=q^{-2m(\alpha+n+1)}\cdot
\frac{(q^{2\alpha+2n+2};q^2)_{m}}{(q^2;q^2)_{m}} \cdot \sum_{l=0}^m
\frac{(q^{2};q^2)_{m}}{(q^2;q^2)_{l}(q^2;q^2)_{m-l}} \cdot (-1)^m\cdot
q^{l(l-2m+1)}\cdot y^{n+1+l+\alpha}=$$
\begin{equation}\label{f}
=q^{-2m(\alpha+n+1)} \cdot \sum_{l=0}^m
\frac{(q^{2\alpha+2n+2};q^2)_{m}}{(q^2;q^2)_{l}(q^2;q^2)_{m-l}} \cdot
(-1)^m\cdot q^{l(l-2m+1)}\cdot y^{n+1+l+\alpha}.
\end{equation}
Hence it remains to compute $\mathscr{F}y^{n+1+l+\alpha}$ for any $l\geq0$.
Let us apply the $q$-spherical transform to the both hand-sides of the
formula $B_{q,\alpha}(f_0)=(q^{2\alpha+2};q^2)_n\cdot y^{\alpha+n+1}$. Due
to Proposition \ref{sy}

$$ \frac{(q^{2+2\alpha};q^2)_\infty\cdot(q^{2n+2+2\alpha};q^2)_\infty}
{(q^{n+2+2\alpha+i\rho};q^2)_\infty\cdot(q^{n+2+2\alpha-i\rho};q^2)_\infty}
\cdot
\mathscr{F}{f}_0=(q^{2+2\alpha};q^2)_n\cdot\mathscr{F}y^{n+1+\alpha},$$ or,
equivalently,

$$\mathscr{F}y^{n+1+\alpha}=\frac{(q^{2};q^2)_n\cdot(q^{2n+2+2\alpha};q^2)
^2_\infty}
{(q^{n+2+2\alpha+i\rho};q^2)_\infty\cdot(q^{n+2+2\alpha-i\rho};q^2)_\infty}.$$
Let us make the substitution $N:=n+1+\alpha$ ($N>n+1$):
$$\mathscr{F}y^{N}=\frac{(q^{2};q^2)_n\cdot(q^{2N};q^2)^2_\infty}
{(q^{2N-n+i\rho};q^2)_\infty\cdot(q^{2N-n-i\rho};q^2)_\infty}.$$ Thus,
$$\mathscr{F}y^{n+1+l+\alpha}=
\frac{(q^{2};q^2)_n\cdot(q^{2n+2+2l+2\alpha};q^2)^2_\infty}
{(q^{n+2+2l+2\alpha+i\rho};q^2)_\infty
\cdot(q^{n+2+2l+2\alpha-i\rho};q^2)_\infty}.$$
The latter equality, together with (\ref{f}) and simple computations, gives
our formula.\hfill$\blacksquare$

\medskip

Recall the notation $p_k(x;a,b,c|q)$ for the continuous dual $q$-Hahn
polynomials \cite[Section 3.3]{KS}. The above proposition implies
\begin{equation*}
\begin{split}
&\qquad \mathscr{F}\sigma^{(\mathrm{r})}(P_m)=
\\
&=\frac{q^{-mn}\cdot(q^{2};q^2)_n\cdot (q^{2n+2+2\alpha};q^2)_\infty\cdot
(q^{2n+2m+2+2\alpha};q^2)_\infty} {(q^{n+2+2\alpha+i\rho};q^2)_\infty
\cdot(q^{n+2+2\alpha-i\rho};q^2)_\infty \cdot(q^{2};q^2)_m}\cdot
p_m(\cos\frac{h\rho}{2};q^{n+2+2\alpha},q^n,q^n|q^2).
\end{split}
\end{equation*}
Let $p_m(\rho):=p_m(\cos\frac{h\rho}{2};q^{n+2+2\alpha},q^n,q^n|q^2).$
Applying Proposition \ref{main} and above observations, we get

$$\frac1{4\pi}\int\limits_0^{\frac{2\pi}{h}}p_m(\rho)p_l(\rho)\frac{d\rho}
{(q^{n+2+2\alpha+i\rho};q^2)^2_\infty\cdot(q^{n+2+2\alpha-i\rho};q^2)^2_\infty
\cdot b_{q,\alpha}(\rho)|c(\rho)|^2}=$$
$$=\frac{q^{2mn}\cdot(1-q^{2n})\cdot(q^{2};q^2)^2_m}
{h\cdot(q^{2};q^2)^2_n\cdot (q^{2n+2+2\alpha};q^2)^2_\infty\cdot
(q^{2n+2m+2+2\alpha};q^2)^2_\infty}\cdot\mathrm{Tr}_q(P_m\cdot P_l),$$ or,
using the explicit formula for $b_{q,\alpha}$,

$$\frac1{4\pi}\int\limits_0^{\frac{2\pi}{h}}p_m(\rho)p_l(\rho)\frac{d\rho}
{(q^{n+2+2\alpha+i\rho};q^2)_\infty(q^{n+2+2\alpha-i\rho};q^2)_\infty
|c(\rho)|^2}=$$$$=\delta_{ml}\frac{(q^{2};q^2)_m(q^{2};q^2)_{m+n-1}}
{h(q^{2};q^2)_{n-1}(q^{2n+2m+2+2\alpha};q^2)^2_\infty}.$$ Recall (Proposition
\ref{in}) that

$$\frac1{|c(\rho)|^2}=\left|\frac{\Gamma_{q^2}^2(\frac{n}{2}+\frac{i\rho}{2})}
{\Gamma_{q^2}(n)\Gamma_{q^2}(i\rho)}\right|^2=(q^{2n};q^2)^2_\infty\cdot
\frac{(q^{2i\rho};q^2)_\infty(q^{-2i\rho};q^2)_\infty}
{(q^{n+i\rho};q^2)^2_\infty(q^{n-i\rho};q^2)^2_\infty}.$$ Finally we have

$$\frac1{4\pi}\int\limits_0^{\frac{2\pi}{h}}p_m(\rho)p_l(\rho)
\frac{(q^{2i\rho};q^2)_\infty(q^{-2i\rho};q^2)_\infty}
{(q^{n+2+2\alpha+i\rho};q^2)_\infty(q^{n+2+2\alpha-i\rho};q^2)_\infty
(q^{n+i\rho};q^2)^2_\infty(q^{n-i\rho};q^2)^2_\infty}d\rho=$$$$
=\delta_{ml}\frac1 {h(q^{2n+2m+2+2\alpha};q^2)^2_\infty(q^{2m+2};q^2)_\infty
(q^{2n+2m};q^2)_\infty}.$$ The latter equality is a particular case of the
orthogonality relations for the entire family of continuous dual $q$-Hahn
polynomials given in \cite[(3.3.2)]{KS}.

\bigskip

\section{Conclusion}\label{7} As it was mentioned in the Introduction, this
research is, above all, a part of the general program of studying $q$-Cartan
domains, which we believe is a promising subject in quantum group theory. To
our knowledge the results above on finding a mathematical setting for the
Toeplitz and covariant calculi, expressing Berezin transform as a function
of the $q$-Laplacian operator, and on computing the covariant symbol of the
projections are all new. Those results, in the classical case, are naturally
related to finding an expansion of the associated $\ast$-product $f\ast_h g$
of covariant symbols as a power series in the Planck constant $h$ (see e. g.
\cite{Oesterle}). Thus, the most natural continuation of the present paper
would be understanding the expansion of the $\ast$-product in this case
since certain computation might be easier than in the classical case, and
generalization of our results to the case of other, more complicated
$q$-Cartan domains, first of all, $q$-matrix balls \cite{SSV1, SSV2}.
 We mention the paper \cite{SSV2} where $q$-weighted
Bergman spaces on $q$-matrix balls are constructed and studied. Using these
results, it is not so difficult to generalize our definition of the
$q$-Berezin transform. However, to obtain more deep results (for instance,
an explicit asymptotic formula for the $q$-Berezin transform) one needs to
develop to some extent harmonic analysis on $q$-matrix balls. This seems to
be a difficult problem itself. By now, harmonic analysis is developed quite
well for compact quantum homogeneous spaces only. Nevertheless, we are fully
confident that all the results we present in the paper admit generalization
to other $q$-Cartan domains.

\end{document}